\documentclass[aop]{imsart}
\usepackage{graphicx}                      

\usepackage{amsfonts}
\usepackage[margin=1.0 in]{geometry}

\RequirePackage{amsthm,amsmath}  

\usepackage{amssymb}

\startlocaldefs
\numberwithin{equation}{section}
\theoremstyle{plain}
\newtheorem{thm}{Theorem}[section]
\newtheorem{lemma}{Lemma}[section]

\newtheorem{cor}[thm]{Corollary}
\endlocaldefs

\long\def\/*#1*/{}  

\usepackage{placeins} 
\usepackage{bm} 
\usepackage{enumitem} 
\usepackage{caption}

\newtheoremstyle{boldhead}
  {\topsep}
  {\topsep}
  {\itshape}
  {}
  {\bfseries}
  {.}
  {.5em}
  {\thmname{#1}\thmnumber{ #2}\thmnote{ (#3)}}

\theoremstyle{boldhead}
\newtheorem{boldthm}[thm]{Theorem}

\begin{document}

\begin{frontmatter}
\title{ Concentration Inequalities for Branching Random Walks }

\begin{aug}
\affiliation{Changqing Liu  \\ ClinTFL Ltd. \\  c.liu@ClinTFL.com}
\end{aug}

\begin{abstract}
Motivated by phase transition problems in CSPs, we prove a more general concentration inequality that retains classical sub-Gaussian tails under a mild global linear-growth condition $|S_n| < Cn$, relaxing the bounded-increment assumption to finite exponential moments and requiring neither independence nor the martingale property.

We further extend it to branching random walks (BRWs), obtaining the first concentration inequality for BRWs.
\end{abstract}   
  

\end{frontmatter}

\section {Concentration Inequalities}
\subsection{Introduction}\label{intro con ine}
\textbf{}Let $(S_i)_{i=1, 2, \ldots}$ be a real-valued random process (a one-dimensional random walk). $S_n$ can be formulated as $S_n=X_1 + X_2 + \ldots + X_n$; $X_i=S_i-S_{i-1}$ are called increments. Under certain conditions, in particular independence, $S_n$ converges to a Gaussian distribution within a $\sqrt{n}$ neighborhood of its mean---the central limit theorem (CLT). For $S_n$ beyond the neighborhood, i.e., $|S_n - \mathbb{E}(S_n)| \gg \sqrt{n}$, Chernoff inequality (\cite{Chernoff 1952} \cite{Hoeffding 1963}),  also referred to as Chernoff bound,
\begin{equation}\label{chernoff bound}
         \Pr(|S_n - \mathbb{E}(S_n)| \geq \lambda) \leq e^{- c\frac{\lambda^2}{n}}
\end{equation}
tells us how unlikely this occurrence is.
Chernoff inequality extends to the setting of bounded martingale difference (e.g., $\mathbb{E}(X_n \,|\, S_{n-1}) = 0$), known as Azuma-Hoeffding inequality  \cite{Azuma}  (see \cite{Chung} and \cite{McDiarmid89} for surveys and references therein). From a practical perspective, the Azuma-Hoeffding inequality and Chernoff inequality give the same tight probability bound: square-exponential decay\footnote{Therefore, the constant `$c$' in the Chernoff bound is immaterial and the same symbol is used for simplicity.}, while the former has wider scope of applicability.
After the classical results of Hoeffding and Azuma(\cite{Hoeffding 1963}, \cite{Azuma}), many refinements and extensions followed (e.g. \cite{Talagrand 1995}, \cite{McDiarmid89}, \cite{Bentkus 2004}, \cite{bentkus 2008}, \cite{Fan & Grama & Q. Liu 2012}, \cite{Kontorovich 2014}). These concentration inequalities have proven very useful and have a wide variety of applications in computer science, combinatorics, information theory (see e.g., \cite{Dubashi&Panconesi 2009}, \cite{Alon2000},  \cite{Raginsky&Sason 2013}, \cite{Motwani&Raghavan 1995}).

For a long time, no substantive progress was made in eliminating the boundedness condition (the trivial idea of truncation gives weak results, though)\footnote {In fact, it was not until 1995 that the first major improvement in Hoeffding inequality appeared \cite{Talagrand 1995}. The reason might be twofold. One reason may be the lack of strong motivation. Another, nifty techniques employed in the pioneering work, and in the subsequent literatures as well, achieving finer results somehow obscured the fundamental mechanism underlying the concentration inequalities.}. In the early 2000s \cite{Lesigne&Volny} relaxed the boundedness to finite exponential moments; i.e. for a $\delta > 0$, $\mathbb{E}(e^{\delta X_i}) $ is uniformly bounded for all $i$. They obtained  
\[
     \Pr(\,|S_n - \mathbb{E}(S_n)| \geq \lambda|\,) \leq \exp\left(-c \Big( \frac{\lambda^2}{n} \Big)^{1/3}\right)
\]
which falls short of a sub-Gaussian bound. Later, motivated by study of the free energy of a directed polymer in a random environment,  \cite{Q_Liu&Watble} gave a Chernoff bound under the hypothesis that, for some $\delta > 0$, $\mathbb{E}(e^{\delta |X_i|} \mid \mathcal{F}_{i-1}) $ is uniformly bounded for all $i$ \big(where $\{\mathcal{F}_i\}$ is a filtration for the martingale difference sequence $\{X_i\}$\big). Other improvements to concentration inequalities in the relaxation of the boundedness condition can be found in \cite{bentkus 2008} \cite{Fan & Grama & Q. Liu 2012} \cite{Kontorovich 2014} \cite{Fan & Grama & Q. Liu 2015}. 

Wormald's differential equation method (DEM) \cite{Wormald Ann. App Prob} is the first general framework that eliminates the the need for independence and martingale assumptions. To our knowledge, no existing work relaxes both independence/martingale and boundedness conditions while achieving sub-Gaussian concentration. 
The first major contribution of this paper is to do so: informally, with finite exponential moments of the increments, the small per-step effect in a random process can be controlled by the Lipschitz continuity of the mean increment with respect to the history, not necessarily by ``almost independence''. Our approach and the differential equation method (DEM) of Wormald (\cite{Wormald Ann. App Prob}, \cite{Wormald in Warsaw}) share similarities in that both dispense with assumptions of independence or martingale differences in the underlying process, instead employing a martingale construction together with smoothness conditions on the expectation. Although in our framework the absence of the boundedness assumption leads to technically distinct arguments, our approach and Wormald's can be viewed as mutually explanatory, each shedding a different light on how concentration arises without independence or martingale-difference assumptions in the underlying process.

\subsection {Notation and first contribution}\label{Notation and first contribution}
Let $U_i$ denote the partial sum of the first $i$ terms of  $S_n = \sum_{j=1}^n X_j$, namely, $U_i = X_1 + X_2 + \cdots + X_i$. We use lower case for variables scaled  by $1/N$, where $N$ is a constant integer (in the context of $K$-SAT, $N$ represents the total number of variables). For instance, $u_i := \sum_{j=1}^i X_j/N$, $x_i := X_i/N$. For the sake of readability, in some places we use $Y_n$ for $S_n$. 

We let $\mathbb{E}_{i}(\,\cdot \,)$ abbreviate $\mathbb{E}(\, \cdot \mid X_1, ..., X_i \, )$. 
To express the stochastic process in calculus notation, we define
\begin{equation}\label{Py/Pt}
    \frac{\partial y} {\partial t} := 
    \frac{\mathbb{E}_i(y_{t}) - \mathbb{E}_i(y_{t-1})}{  1/N}
     =  \mathbb{E}_i(Y_{t}) - \mathbb{E}_i(Y_{t-1}) = E_i(X_{t}) 
     , \quad (\text{for $ t > i$})
\end{equation}
where $Y_t= X_1 + X_2 + \cdots + X_t $. In the martingale-difference case, by the law of iterated expectations, we have  $E_i(X_t) = 0$ for $ t > i$, and thus
\[ 
    \displaystyle \frac { \partial y } {\partial t} = 0
\]
In the case of independence, $\mathbb{E}_i(X_t) = \mathbb{E}(X_t)$, i.e. 
\[
   \displaystyle \frac { \partial y } {\partial t}= \text{constant (with respect of $u_i$)}
\] 
Thus in  both cases, $\partial y / \partial t$ is independent of $u_i$, and we may write
\begin{equation}\label{vanishing mixed derivative}
              \frac { \partial } {\partial u_i} \left( \frac { \partial y} {\partial t} \right)
              =
              \frac { \partial^2 y } {\partial u_i \, \partial t} =0
\end{equation} 
Here, for $t > i$, we define $\mathbb{E}(X_t \mid u_i) := \mathbb{E}\big[X_t \mid \sigma(u_i)\big]$ which is a function of $u_i$ 
\vspace{0.1\baselineskip} 
by Doob–Dynkin lemma.

While classical concentration inequalities require the condition (\ref{vanishing mixed derivative}),  one of the primary contributions of this work is the relaxation of this classical constraint to the condition:
\[  \left| \frac { \partial^2 y } {\partial u_i \, \partial t} \right| \le L, \quad i=1, 2, \ldots, n
\]
i.e., 
$
\displaystyle \frac { \partial^2 y } {\partial u_i \, \partial t} 
$ 
is required only to be uniformly bounded, thereby allowing the drift $\displaystyle\frac { \partial y } {\partial t}$ to vary with $u_i$. This extends concentration results to general dependent processes with variable drift. 
For instance, the generation of random 3-SAT instances can be modelled as a branching random walk where, at each time step, a set of random clauses is applied to the current formula, producing $O(N^3)$ progeny formulae from the current state; for details, see later in this paper. In this process, certain properties of the formula, such as the number of frozen variables, are history-dependent.

As a side benefit, this translation of classical probabilistic conditions from expectation-based to derivative-based formulations provides a calculus-based perspective on stochastic concentration phenomena.

%
%

%
%

\section{Extending the Hoeffding-Azuma inequality}\label{Extending the Hoeffding-Azuma inequality}
As a pedagogical warm-up, we prove a ``non-sharp'' version of Azuma inequality without assuming bounded increments. 
To begin with, we prove a simple inequality drawing on the method introduced in \cite{Q_Liu&Watble} and \cite{Watbled2012}. 
\begin{lemma}\label{Lemma: Q_Liu&Watble} 
If $\mathbb{E}(X)=0$, then for all $\delta > 0$ and $|\, t \,|/\delta \in [0,1]$
\[
     \mathbb{E} \left[ e^{tX} \right] 
                                  \le \exp \left (t^2/\delta^2 \cdot \mathbb{E} \big[ e^{\delta |X|} \, \big] \right)
\]
\end{lemma}
In particular, when $\delta = 1$,
\begin{equation}\label{E<t^exp} 
\mathbb{E}\left[ e^{tX} \right ] \leq \exp{ \left( t^2\mathbb{E} \big[ e^{|X|} \big]  \right) }   
\end{equation}
{\allowdisplaybreaks
\begin{proof}
   \begin{align*}
   \mathbb{E}(e^{tX}) &= 1 + \mathbb{E}\left[ tX \right]  
                    +  \mathbb{E} \left[ \frac{t^2 X^2}{2!} + \frac{t^3 X^3}{3!} + \cdots \right] 
                                                            \\[5pt]
      & =  1 +  \mathbb{E} \left[ \frac{(t/\delta)^2 \cdot (\delta X)^2}{2!} +
                        \frac{(t/\delta)^3 \cdot (\delta X)^3}{3!}  +  \cdots  \ \right]
                                                                                 \\[5pt]
      &\leq  1 +  (t/\delta)^2\mathbb{E} \left[ \frac{|\delta X|^2}{2!} + \frac{|\delta X|^3}{3!} + \cdots \ \right]
                                   \quad    \text {\big(employing $ |\, t \,|/\delta \leq 1 $\big)}
                                                                                 \\[5pt]
      &\leq  1 +  (t/\delta)^2\mathbb{E} \left [ e^{\delta |X|} \right]         
                                                                                 \\[5pt]      
      &\le \exp \left (t^2/\delta^2 \cdot \mathbb{E} \big[ e^{\delta |X|} \, \big] \right) 
   \end{align*}
The proof is concluded.
\end{proof}
}%

\begin{thm}\label{Azuma}
Let $S_n = X_1 + X_2 + \cdots + X_n$, and $\{X_i\}_{i=1}^n$ be a sequence of martingale differences {\normalfont(thus $\mathbb{E} S_n = 0$)} with $\mathbb{E}[ e^{|X_i|} ] \le K$ for all $i = 1, \dots, n$. 

Then, for every $\lambda > 0$, 
\begin{equation}\label{Pr < min}
       \Pr(|S_n - \mathbb{E} S_n| \ge \lambda) \le 2\exp 
                                  \left(
                                           - \min  \left \{\frac{\lambda^2}{4Kn}, \frac{\lambda}{2} \right \}
                                  \right)    
\end{equation}            
\end{thm}
\noindent This is equivalent to 
\[  
                \Pr(|S_n - \mathbb{E} S_n| \ge \lambda) \le
\begin{cases}
              2\exp  \left(  - \displaystyle \frac{1}{4K}\frac{\lambda^2}{n}  \right)      & \lambda \le 2nK    \\[18pt]                   
      \displaystyle    2\exp  \left(  - \displaystyle \frac{\lambda}{2}  \right)  & \lambda > 2nK   
       
\end{cases}
\] 
since
\[
              \lambda \le 2nK  \iff  \frac{\lambda^2}{4Kn} \le \frac{\lambda}{2}
\]                 
\/**********************************************************************************************************************
To prove the inequality, Azuma (\cite{Azuma}) employed the following inequalities, for $|X| \le c $,
\begin{equation}\label{cosh}
      e^{tX} \le \frac {e^{tc} - e^{-tc}}{2c}X + \frac {e^{tc} + e^{-tc}}{2} \quad \text{(due to convexity of $e^{tX}$)}
\end{equation}
and
\begin{equation}\label{et_e2t}
         \frac {e^{tc} + e^{-tc}}{2} \le e^{t^2c^2/2}  \quad \text{($\cosh(x) < e^{x^2/2}    $)}
\end{equation}
to get sub-gaussian condition 
\[
     \mathbb{E}(e^{tX}) \le e^{t^2c^2/2}  ~~~~~~ \mbox {for}~ \mathbb{E}(X) = 0
\]
as adopted by textbooks. Sub-gaussian condition is a key ingredient in the proof of concentration inequalities, as this  ultimately leads to a bound of the form 
\[   
   \exp({\sim n t^2 - t\lambda}),\quad 
\]
for $\lambda \sim n^{(1+\varepsilon)/2}$ with $0 < \varepsilon < 1$, 
which 
can be made exponentially small---of order $\exp(-\Th eta(n^\varepsilon))$---by choosing sufficiently small $t$. Intuitively, provided $\lambda$ is not too small, the exponential term can be made exponentially small by choosing small $t$, as $t^2$ is of higher order than $t$ and thus $nt^2$ is vanishing faster than $t\lambda$ when $t$ is approaching $0$. 
\begin{proof}[\textnormal{\textsc{Proof} (\text{Extended Azuma inequality}})]
******************************************************************************************************************/
%

\begin{cor}[Extended Azuma inequality]\label{coro Azuma} Suppose, in addition, that $|S_n| < Cn$ for some constant $C > 0$ {\normalfont (as in $K$-SAT, where $S_n$ is the number of backbones; indeed, a formula with $n$ $K$-clauses contains at most $Kn$ backbones)}. Then, for every $\lambda > 0$,
\begin{equation}
       \Pr(|S_n - \mathbb{E} S_n| \ge \lambda) \le 2\exp 
                                  \left(
                                           -  c\frac{\lambda^2}{n}  
                                  \right)                                   
\end{equation}  
where 
\[
       c = \min\left\{\frac{1}{4K}, \frac{1}{2C} \right \}
\]     
\end{cor}
Henceforth, the condition $|S_n| < Cn$ is assumed throughout.  

\begin{proof} [Proof of Theorem \ref{Azuma}]   
We only consider ${ \Pr }(S_n - \mathbb{E} S_n  \ge \lambda)$, as the proof of the reverse ${ \Pr }(S_n - \mathbb{E} S_n  \le -\lambda)$ is exactly the same. First, let $\lambda \le 2nK$. By Markov's inequality, inductive application of Lemma \ref{Lemma: Q_Liu&Watble} and optimization of $t$, we have
{\allowdisplaybreaks
\begin{align*}
      { \Pr }(S_n - \mathbb{E} S_n  \ge \lambda)  &\le e^{-t\lambda} \cdot \mathbb{E} \left[ e^{t(S_n - \mathbb{E}S_n)} \right]  \quad  \text{\big(Note $\mathbb{E} S_n=0$\big)} 
      \\[6pt]       
       &\le
       e^{-t\lambda} \cdot \mathbb{E} \Big[
                                           e^{tS_{n-1}} \cdot \mathbb{E} (
                                           e^{tX_n} \mid X_1, \ldots, X_{n-1}
                                      )
                            \Big]    
       \\[6pt]
       &\le e^{-t\lambda} \cdot \mathbb{E}\Big[
                        e^{tS_{n-1}} \cdot \exp{ \Big( t^2\, \mathbb{E}  \, [\, e^{|X_n|} \mid X_1, \ldots, X_{n-1} \, 
                                                                          ]   \Big) }
                      \Big] \tag{By (\ref{E<t^exp})  of  Lemma \ref{Lemma: Q_Liu&Watble}} 
                                    \\[8pt] 
       &\le e^{-t\lambda} \cdot \mathbb{E}\Big[
                        e^{tS_{n-1}} \cdot e^{t^2 K}
                      \Big]   \tag{$\mathbb{E}[e^{|X_i|}] \le K$}    
      \\
      & \tag{Induction over $n$} \le   \exp (-t\lambda \, + \, nt^2 K) 
      \\
      & = \exp{ \left ( - \frac{\lambda ^ 2}{2n} \cdot \frac{1}{2K}  \right )  } \tag{Optimal $t = \frac{\lambda}{2nK}$}  
      \\  &\hspace{0em}  
           \text{
                  \big(cf. Azuma's  $ \exp { \left(\displaystyle - \frac{\lambda ^ 2}{2n} \cdot \frac{1}{c^2}\right)} $ 
                   \!\!, assuming  $|X_i| \le c $ 
                  \big)
                } 
\end{align*}
}%
Note that $\lambda \le 2nK$ ensures $t \le 1$, as required in Lemma \ref{Lemma: Q_Liu&Watble}. 

\vspace{.5\baselineskip}
For $\lambda > 2nK$ $\big($and hence $ -\lambda + nK < -\lambda/2$$\big)$, Markov's inequality together with $\mathbb{E}[e^{|X_i|}] \le K$  and the tower property yields
{\allowdisplaybreaks
\begin{align*}
    \Pr (S_n - \mathbb{E} S_n \ge \lambda ) & \le   e^{-\lambda} \mathbb{E} e^{S_n}     \le   e^{-\lambda}   K^n      
                                               =   \exp(-\lambda + n\ln K) 
     \\[5pt]
     &\le    \exp(-\lambda + n K)       
     \\[5pt]
     & \le \exp(-\lambda/2)  \qquad \left(\text{since $-\lambda + nK < -\lambda/2$} \right)   
\end{align*}
Combining this with the result for $\lambda \le 2nK$ completes the proof. 
}%
\end{proof} 
  
\begin{proof}[Proof of Corollary \ref{coro Azuma}]   
Bear in mind that $ \mathbb{E} S_n = 0 $. For $\lambda > Cn $, the proof is trivial, since $|S_n| < Cn$ and hence the left-hand side of the concentration inequality is zero. For $ 0 < \lambda \le Cn $, from
\[
         c \le \frac{1}{4K} \text{ and } c \le \frac{1}{2C}
\]
we have 
\[
         c\frac{\lambda^2}{n} \le \frac{\lambda^2}{4Kn}
\]
and
\[
   c\frac{\lambda^2}{n} \le \frac{\lambda^2}{2Cn} \le \frac{\lambda}{2} 
\] 
Therefore 
\[ 
        c\frac{\lambda^2}{n} \le \min \left \{   
                                                 \frac{\lambda^2}{4Kn}, \frac{\lambda}{2} 
                                      \right \}
\]
Hence 
\[
                            \exp\left(    - \min \left \{   
                                                 \frac{\lambda^2}{4Kn}, \frac{\lambda}{2} 
                                      \right \}
                                \right)      
                       \le   \exp \left(
                                      -c\frac{\lambda^2}{n} 
                                  \right)    
\]                                      
By Theorem (\ref{Azuma}), the desired follows. 
\end{proof}

The preceding exposition foregrounds the fact that the condition $X_i < \infty $ can be relaxed to $\mathbb{E}(e^{|X_i|}) < \infty $ without compromising the strength of the results in the sense that from a practical point of view, the difference between the constant factors in the exponent, i.e., $\displaystyle \exp{\left(- \frac{\lambda ^ 2}{2n}\cdot \frac1 {2K}\right)}$ versus $\displaystyle \exp{\left(- \frac{\lambda ^ 2}{2n}\cdot \frac 1 {c^2} \right)}$, is immaterial---both are square-exponential decay. In addition, it is demonstrated that the sub-Gaussian tail bounds can be derived directly from the finite exponential moment without relying on ``nifty tricks'' (such as Hoeffding's Lemma); particularly when bounded increments are considered (i.e., $|X_i| \le c$), we see that the sub-Gaussian bound is obtained straightforwardly this way, taking the form $\displaystyle \exp{\left(- \frac{\lambda ^ 2}{2n}\cdot \frac1 {2e^c}\right)}$.

Next, we identify weaker conditions that ensure concentration inequalities, beyond independence, martingale differences, and boundedness. We begin by defining Doob's (or McDiarmid's) martingale, 
\begin{align*}\label{Doob's} 
  \Delta_i &= \mathbb{E}_{i}(X_i + \cdots + X_n) - \mathbb{E}_{i-1}(X_i + \cdots + X_n)
\end{align*}
where $\mathbb{E}_i(\,\cdot\,) := \mathbb{E} (\, \cdot \mid X_1, ..., X_i) $,  and similarly $\mathbb{E}_{i-1}(\,\cdot\,) := \mathbb{E} (\, \cdot \mid X_1, ..., X_{i-1}) $.
It is easy to check
\[       S_n - \mathbb{E} S_n = \sum_{i=1}^n \Delta_i \,
\begingroup
\renewcommand{\thefootnote}{$\S$}
\footnote{  $ \mathbb{E}_n[S_n] - \mathbb{E}_0[S_n] = S_n - \mathbb{E}S_n $    
} 
     \addtocounter{footnote}{-1}     
\endgroup%
                    \, , \quad \mathbb{E} _{i-1}(\Delta_i) = 0
\] 
i.e. $\Delta_i$ is a martingale difference. 
Note 
\begin{align} \label{di}  
      \Delta_i &=  X_i- \mathbb{E}_{i-1}X_i+\sum_{\ell = i+1}^n\big( \mathbb{E}_i X_\ell - \mathbb{E}_{i-1}X_\ell \big) 
\end{align}

%
%
%

%
%
%
%
 
\begin{lemma}[L-Lipschitz]\label{Lip_f} If 
\begin{equation}\label{L|X|/N}
  \left| \,  
          \mathbb{E}_{i-1}(X_\ell \mid X_i) - \mathbb{E}_{i-1}(X_\ell \mid X_i')  
        \, 
  \right | 
              \le L \displaystyle \frac{|X_i - X_i'|}{N} 
              ,\quad \text{for $i < \ell $}  
\end{equation} 
then
\begin{equation}\label{E_LXi} 
                    |\mathbb{E}_i(X_\ell)-\mathbb{E}_{i-1}(X_\ell) | 
                    \le L \frac{|X_i|}{N} + L \frac{\mathbb{E}_{i-1}(\,|X_i|\,)}{N}
                    ,\quad \text{for $ i < \ell $}  
\end{equation}
 
\end{lemma}   

\begin{proof}[\textnormal{\textsc{Proof} (\textsc{Lemma}} \ref{Lip_f})]
   
Let $X_i'$ be a conditionally independent copy of $X_i$ given $(X_1,\ldots,X_{i-1})$; that is,
\[
      X_i'\sim \mathcal{L}(X_i\mid X_1, ..., X_{i-1})
\]
Define
\[
     \mathbb{E}_{X_i'}(\,\cdot\,) \equiv \mathbb{E}_{X_i'}(\, \cdot \mid X_1, ..., X_{i-1})      
\]  
where $\mathbb{E}_{X_i'}(\cdot)$ averages over $X_i'$.  
With this notation and the tower property of conditional expectation, we have, 
{\allowdisplaybreaks
\begin{align*}
      \bigl|\mathbb{E}_i \left( X_\ell \right) - \mathbb{E}_{i-1} \left( X_\ell \right)\bigr| 
      &= \Bigl |\, \mathbb{E}_{i-1}\left( X_\ell \mid X_i \right)  
                         -  
                   \mathbb{E}_{X_i'}\!\left[\,   
                                                 \mathbb{E}_{i-1}\!\left(\, X_\ell \mid X_i'\right)   
                                      \right]
         \Bigl | 
      \\[5pt] &
       = \Bigl |\, \mathbb{E}_{X_i'} \big[\,
                   \mathbb{E}_{i-1} \left( X_\ell \mid X_i \right)  -  \mathbb{E}_{i-1}\!\left(\, X_\ell \mid X_i'\right)  
                                     \big]
        \, \Bigl | 
      \\[5pt] &
       \le \mathbb{E}_{X_i'} \Big[\,
       \bigl |\,  
                 \mathbb{E}_{i-1} \left( X_\ell \mid X_i \right )  -  \mathbb{E}_{i-1}\!\left(\, X_\ell \mid X_i'\right)  
         \bigl |
                             \, \Big]
      \\[5pt] &                             
       \le \mathbb{E}_{X_i'}\! \left[\,   
                                       L \frac{|X_i - X_i'|}{N}  
                             \, \right] 
         \qquad (\text{by \eqref{L|X|/N}})
      \\[5pt] & 
       \le                     
       \frac L N   |X_i| +  \frac L N   \mathbb{E}_{X_i'} \left(\, |X_i'| \,\right) 
      \\[5pt] & 
      =  L \frac {|X_i|} N +  L \frac {\mathbb{E}_{i-1} \left(\, |X_i| \,\right)   } N      
\end{align*}      
}%
The claim is true.
\end{proof}


\begin{lemma}\label{con_extension}  
Let $A$ denote $1+L\frac{M}{N}$, where $M=O(N)$. If the Lipschitz condition is satisfied, i.e. 
\begin{equation}\label{Lipschitz}
  \left| \,  
          \mathbb{E}_{i-1}(X_\ell \mid X_i) - \mathbb{E}_{i-1}(X_\ell \mid X_i')  
        \, 
  \right | 
              \le L \displaystyle \frac{|X_i - X_i'|}{N} 
              ,\quad \text{for $i < \ell $}  
\end{equation} 
and 
\[
       \mathbb{E} [\, e^{\delta A|X_i|}\mid u_{i-1} \, ] \le K  \quad \textnormal{(Finite MGF)} 
\]       
for a constant $\delta > 0$ , $i=1, 2, \ldots , n$, 
then for $n \le M$\renewcommand{\thefootnote}{$\dagger$} \footnotemark 
\begingroup
                   \renewcommand{\thefootnote}{$\dagger$}
                  \footnotetext{Actually, the lemma holds for $n = \Theta(M)$; we assume $n \le M$ here to simplify the argument. 
                               }
                  \addtocounter{footnote}{-1}%
\endgroup 
\begin{equation}\label{Extension for RW}
       \Pr(|S_n - \mathbb{E} S_n| \ge \lambda) \le 2\exp 
               \left(
                     - \min  \left \{
                                       \frac{\delta^2}{4K^2} \frac{\lambda^2}{n}, \frac{\delta \lambda}{2} 
                            \right \}
              \right)  
\/************************                                      
              \renewcommand{\thefootnote}{$\S$} \footnotemark          
              \begingroup
                   \renewcommand{\thefootnote}{$\S$}
                  \footnotetext{Generally, $c(\delta)$ admits a maximum, denoted by $c(\delta^\star)$; heuristically, $c(\delta)$ is close to zero for both very small and very large $\delta$. Hence 
                     $c(\delta^\star) \ge c(1)  \ge  \displaystyle {1}/{(4{K_M}^2) }, 
                                                       $
                   where $ \displaystyle K_{M} := \max_{i}  \mathbb{E} \left[\, e^{A|X_i|}\mid u_{i-1} \, \right]$.
                               }
                  \addtocounter{footnote}{-1}%
              \endgroup  
**************************/              
\end{equation}
\end{lemma}

\begin{cor}\label{coro con_extension} Suppose, in addition, that $|S_n| < Cn$ for some constant $C > 0$ {\normalfont(as in $K$-SAT discussed in this paper, and many other CSPs, including all cases where Wormald's DEM applies)}. 
%
Then, for every $\lambda > 0$,  
\begin{equation}\label{con ineq sub Gaussian}
       \Pr(|S_n - \mathbb{E} S_n| \ge \lambda) \le 2\exp 
                                  \left(
                                           -  c\frac{\lambda^2}{n}  
                                  \right)                                   
\end{equation}  
where 
\[   
       c = \min \left\{\frac{\delta^2}{4K^2}, \frac{\delta}{4C}  \right \}       
\]     
\end{cor}
Throughout the remainder of the paper, by ``concentration inequality'' we mean the sub-Gaussian tail bound \eqref{con ineq sub Gaussian}.

\textbf{Remarks}. 
\begin{itemize}

\item The finite MGF assumption is equivalent to the Bernstein-type moment condition $\mathbb{E}|X|^k < k!\, K^k$, up to the choice of constants, but more user-friendly than the latter; in practice, Bernstein-type bounds are often derived from the existence of an exponential moment.

\item Condition (\ref{Lipschitz}) in the Lemma \ref{con_extension} is equivalent to the following Lipschitz continuity
\begin{equation}\label{Lipschitz continuity} 
       \left| \mathbb{E}(X_\ell \mid u_i) - \mathbb{E}(X_\ell \mid u_i') \right| \le L\left| u_i - u_i' \right| 
                         \quad \mbox{for $i < l$} 
\end{equation}
where 
\[ 				u_i  = \displaystyle \frac{1}{N}  (X_1 + \cdots + X_{i-1} + X_i), \quad
          u_i' = \displaystyle \frac{1}{N}  (X_1 + \cdots + X_{i-1} + X_i').
\]

Informally, the Lipschitz condition (\ref{Lipschitz}) states that the drift is Lipschitz continuous with respect to history, including the initial state $S_0$. We will later show this holds for smooth drift functions. 

\item
   Independence and martingale difference assumptions both imply the L-Lipschitz condition (\ref{Lipschitz}). That is, the L-Lipschitz condition is weaker than those assumptions.

\item
\noindent Poisson$(d)$, Binomial$(N, d/N)$, and the variable occurrence distribution in SAT/COL models, for example, all satisfy the finite exponential moment condition, as does the Gaussian distribution.

\item 
  For $n = O(N)$,  the L-Lipschitz condtion may imply
                               \[ \displaystyle
                                  \left | \sum_{\ell=i+1}^n \mathbb{E}_{i-1}(X_\ell \mid X_i) 
                                        - \sum_{\ell=i+1}^n \mathbb{E}_{i-1}(X_\ell \mid X_i')   
                                   \right| 
                                \le  L \displaystyle {|X_i - X_i'|} 
                               \]
             When the increment $X_i$ (and $X_i'$) are bounded, this reduces to McDiarmid's bounded difference condition \cite{McDiarmid89}. 
     
\item For sufficiently small $\omega$, $\Delta_i = \mathbb{E}_{i}(X_i + \cdots + X_\omega) - \mathbb{E}_{i-1}(X_i + \cdots + X_\omega)$ is bounded if each $X_i$ is bounded. In this case, by the Azuma–Hoeffding inequality, $ \sum_{i=1}^\omega \Delta_i$ satisfies a sub-Gaussian concentration bound. The whole process of $n$ steps can be segmented into a number of small pieces of size $\omega$\,; each is concentrated around its expectation, and thus the entire process exhibits concentration. A ``piecewise concentration'' argument of this kind is used in Wormald's DEM framework (\cite{Wormald Ann. App Prob},
\cite{Wormald in Warsaw}).

\end{itemize}

\begin{proof}[Proof of Lemma \ref{con_extension}]   
First, consider the small $\lambda$ regime, i.e. 
\[ 
          \lambda \le   2nK^2/\delta 
\] 
By definition (\ref{di}), we have
{\allowdisplaybreaks
\begin{align*}
              &\mathbb{E}_{i-1}  \left[ e^{  \delta  |\Delta_i|   }  \right]
               = \mathbb{E}_{i-1} \exp \left( \delta \, 
          \Big |  
                  X_i- \mathbb{E}_{i-1}X_i \,+ \sum_{\ell = i+1}^n\big( \mathbb{E}_i X_\ell - \mathbb{E}_{i-1}X_\ell \big) 
        \, \Big| 
                                         \right) 
                                       \\[6pt]
                      &\le \mathbb{E}_{i-1}
             \exp \left( 
                         \delta | X_i |  \,+\, \delta | \mathbb{E}_{i-1}(X_i) | 
                               \,+\,  \delta \sum_{\ell=i+1}^n \left|\, \mathbb{E}_i X_\ell -\mathbb{E}_{i-1} X_\ell  \,\right|
                  \right)        
                                       \\[6pt]
\tag{Lemma \ref{Lip_f}}                      & \le e^{\mathbb{E}_{i-1}(\delta |X_i|)} \, \mathbb{E}_{i-1} 
                            \exp   \left (
                                         \delta | X_i |  \,+\,
                                         \delta \sum_{\ell=i+1}^n L \frac{|X_i|}{N}  
                                                         \,+\,
                                         \delta \sum_{\ell=i+1}^n L \frac{ \mathbb{E}_{i-1} |X_i|}{N} 
                                     \right )    
                                       \\[6pt] 
                      & \le e^{\mathbb{E}_{i-1}(\delta |X_i|)}   \, \mathbb{E}_{i-1} 
                              \exp \left(
                                         \delta |X_i| \,+\,
                                         \delta |X_i| L\frac{M}{N} 
                                                      \,+\,
                                         \delta \, \mathbb{E}_{i-1}(\,|X_i|\,)\cdot L \frac{M}{N}
                                   \right) 
                                       \\[6pt]
                      & = \exp \left( \delta \mathbb{E}_{i-1} |X_i| \, + \, \delta L\frac{M}{N} \mathbb{E}_{i-1}|X_i| \right)
                          \cdot
                             \mathbb{E}_{i-1}  
                                                \exp \left( \delta |X_i| \Big( 1 + L\frac{M}{N} \Big) 
                                                    \right)   
                                       \\[6pt]
                      & = \exp \big(\delta A \mathbb{E}_{i-1} |X_i| \big) 
                             \cdot \mathbb{E}_{i-1} \left[  e^{\delta A |X_i| }  \right] \le K^2   
                                      \\  &\hspace{10em}  \text{(Jensen's inequality and the hypothesis)}  
\end{align*}
}%
Thus, for every $i$  
\begin{equation}\label{<= K^2}
                \mathbb{E}_{i-1}  \left[ e^{  \delta  |\Delta_i|   }  \right]       \le K^2   
\end{equation} 
Applying the exponential Markov inequality $\Pr[X \ge \lambda] \le e^{-t\lambda}\, \mathbb{E} (e^{tX})$ and noting
                 $  \mathbb{E}_{i-1}(\Delta_i) = 0   $, 
we have
{\allowdisplaybreaks
\begin{align*}          
      \Pr \big( S_n - \mathbb{E} S_n   \ge \lambda \,\big)   
    &\le  e^{-t\lambda}\,
    \mathbb{E}\!\left[ e^{\,t(S_n - \mathbb{E}S_n)} \right]  
    \\[0.5em] 
    &=  e^{-t\lambda}\,
    \mathbb{E}\!\left[
       e^{t(\Delta_1 + \cdots + \Delta_{n-1})}\;
       \mathbb{E}_{n-1}  \left( e^{t \Delta_n} \right)
    \right]                                                    \\[0.5em]  
    &\le
    e^{-t\lambda}\,
    \mathbb{E}\!\left[
       e^{t(\Delta_1 + \cdots + \Delta_{n-1})}\;
           \exp \Big( t^2/\delta^2 \cdot \mathbb{E}_{n-1} \big[\, e^{\delta |\Delta_n|} \, \big]\, \Big)  
                \right]
    \quad\text{(by Lemma \ref{Lemma: Q_Liu&Watble})}           \\[0.5em] 
    &\le
    e^{-t\lambda}\, e^{\,t^2 K^2/\delta^2}\;
    \mathbb{E}\!\left[
       e^{t(\Delta_1 + \cdots + \Delta_{n-2})}\;
       \mathbb{E}_{n-2}\!\left(e^{t \Delta_{n-1}} \right)
    \right]
    \quad \text{(by (\ref{<= K^2}))}
    \\
    &\quad \vdots
    \\
    &\le
    e^{-t\lambda}\, e^{\,n t^2 K^2/\delta^2}. 
\end{align*} 
}%
The optimal $t$ gives
\[ {\Pr} (S_n - \mathbb{E} S_n \ge \lambda)
            \le \exp \left( - \frac{\delta^2}{4K^2} \frac{\lambda^2}{n} \right) 
\]   
The value of the optimal $t$ is  
\[\displaystyle
    t = \frac{\lambda}{2nK^2/\delta^2}  
\] 
Thus
\[
      \frac{t}{\delta} = \frac{\lambda \delta}{2nK^2}
\]
Hence, in the small $\lambda$ regime, i.e., 
\[ 
            \lambda \le   \frac{2nK^2}{\delta }    
\]
the optimal $t$ satisfies the hypothesis $t/\delta \le 1$ of Lemma \ref{Lemma: Q_Liu&Watble}. 

For the other side of the concentration inequality, note that Lemma \ref{Lemma: Q_Liu&Watble} holds for negative $t$ as well. The proof is then identical to the above.  
\\

Next, in the large $\lambda$ regime, 
\[
    \lambda >  \frac{2nK^2}{\delta }  
                 \text{, \  which implies \ } -\delta \lambda + n {K^2}   < - \frac{\delta\lambda}{2}, 
\] 
Note that \eqref{<= K^2} 
\[
                \mathbb{E}_{i-1}  \left[ e^{  \delta  |\Delta_i|   }  \right]       \le K^2   
\]
holds independently of $\lambda$. 
The exponential Markov inequality and the law of total expectation give
{\allowdisplaybreaks
\begin{align*}          
      \Pr \big( S_n - \mathbb{E} S_n   \ge \lambda \,\big)   
    &\le  e^{-\delta \lambda}\,
    \mathbb{E}\!\left[ e^{\,\delta (S_n - \mathbb{E}S_n)} \right]  
     =  e^{-\delta \lambda}\,
    \mathbb{E}\!\left[
       e^{\delta (\Delta_1 + \cdots + \Delta_{n-1})}\;
       \mathbb{E}_{n-1}  \big( e^{\delta  \Delta_n} \big)  
    \right]   
    \\[0.5em]   
    & \le  e^{-\delta \lambda}\,
    \mathbb{E}\!\left[
                             e^{\delta (\Delta_1 + \cdots + \Delta_{n-1})}\;  K^2   
               \right] \tag{By \eqref{<= K^2}}   
    \\[0.5em]   
    & \le  e^{-\delta \lambda}\,           (K^2)^n   \tag{Induction on $n$} 
      =  \exp( -\delta \lambda\, + n\ln K^2 )         
    \\[0.5em]   
    &  <  \exp( -\delta \lambda\, + nK^2 )   
    \\[0.5em]   
    &  <  \exp \left( - \frac{\delta\lambda}{2}  \right)        
\end{align*} 
}%
The argument for the opposite side is identical. 

Combining the estimate of the small $\lambda$ case, we have 
\begin{equation} 
       \Pr(|S_n - \mathbb{E} S_n| \ge \lambda) \le 2\exp 
               \left(
                     - \min  \left \{
                                       \frac{\delta^2}{4K^2} \frac{\lambda^2}{n},\; \frac{\delta \lambda}{2} 
                            \right \}
              \right)    
\end{equation} 
completing the proof. 
\end{proof}

\begin{proof} [Proof of Corollary \ref{coro con_extension}]   

It is easy to check that the global bound $|S_n| < Cn$ implies  $|S_n - \mathbb{E} S_n| < 2Cn$. For $\lambda > 2Cn $, the proof is trivial, since the left-hand side of the concentration inequality is zero. We now consider the case $ 0 < \lambda  \le 2Cn $. Recall that 
\[   
       c = \min \left\{\frac{\delta^2}{4K^2},\, \frac{\delta}{4C}  \right \}       
\]     
Thus,
\[
         c \le \frac{\delta^2}{4K^2} \text{ \ and \ } c \le \frac{\delta}{4C}
\]
We have 
\[
         c\frac{\lambda^2}{n} \le \frac{\delta^2}{4K^2} \frac{\lambda^2}{n}
\]
and, using $ 0 < \lambda  \le 2Cn $, 
\[
   c\frac{\lambda^2}{n} \le  \frac{\delta}{4C} \frac{\lambda^2}{n} \le \frac{\delta \lambda}{2} 
\] 
Therefore 
\[ 
        c\frac{\lambda^2}{n} \le \min \left \{   
                                               \frac{\delta^2}{4K^2} \frac{\lambda^2}{n}, \frac{\delta \lambda}{2} 
                                      \right \}
\]
and hence
\[
               \exp\left(  - \min \left \{   
                                             \frac{\delta^2}{4K^2} \frac{\lambda^2}{n}, \frac{\delta \lambda}{2} 
                                \right \}
                 \right)      
                       \le   \exp \left(
                                      -c\frac{\lambda^2}{n} 
                                  \right)    
\]                                      
By Lemma   (\ref{con_extension}), the desired follows. 

\end{proof}

\textbf{Remarks}. 

\begin{itemize} 

\item  The Lipschitz condition (\ref{Lipschitz continuity}) may be interpreted as Lipschitz continuous
dependence of the trajectory on the history; e.g., small perturbations in the initial position do not cause drastic changes in the trajectory.

\item  In the case of independence, $\mathbb{E}_i[ X_\ell ] = \mathbb{E}[X_\ell]$, and hence $\mathbb{E}_i[X_\ell]-\mathbb{E}_{i-1}[X_\ell]$ = 0. In the case of martingale difference sequence, $\mathbb{E}_i[X_\ell]=0$ for $\ell>i$. In both cases,
 \[
      \sum_{\ell > i} \big( \mathbb{E}_i X_\ell - \mathbb{E}_{i-1}X_\ell \big) = 0
 \] 
and therefore $\Delta_i$ (see, (\ref{di})) reduces to $X_i - \mathbb{E}X_i$ (independence) and to $X_i$ (martingale difference), respectively. In these cases, the condition of bounded difference or the existence of exponential moment alone suffices to imply the concentration inequalities. 

\item  Recall the notation introduced in the Subsection \ref{Notation and first contribution}. The condition (\ref{Lipschitz continuity}) is 
\begin{equation} \label{cond: Lipschitz} 
       \left|    \frac { \partial } {\partial u_i} \left( \frac { \partial y } {\partial t} \right)  \right| \le L .  
\end{equation} 

\vspace{0.5em} 
\item Two special cases of (\ref{Lipschitz continuity}) are one of independence where $\displaystyle \frac { \partial y } {\partial t}$ is constant (w.r.t. $u_i$), and one of martingale difference where  $ \displaystyle \frac { \partial y } {\partial t} = 0$.

\end{itemize}
 
The following theorem provides readily verifiable conditions for concentration inequality in general settings where
$$ 
        \frac { \partial } {\partial u_i} \left( \frac { \partial y } {\partial t} \right) \ne 0 
$$


\FloatBarrier 
\begin{figure}[ht]    
    \centering       
    \includegraphics[height=13.5em, width=18em, trim=10 4 35 14, clip]{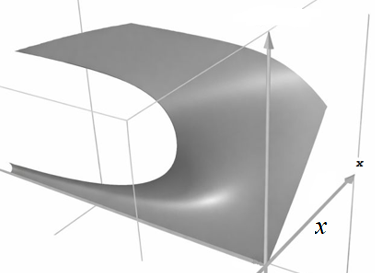}
    \caption{$u(x,z)$, 3-SAT \cite{Liu K-SAT/q-COL}}
    \label{fig: 3-SAT}  
\end{figure}

\vspace{1\baselineskip}
Let $\mu_n(u_n)$ denote $\mathbb{E} (X_{n+1} \mid u_n)$\!\!\!
\begingroup
\renewcommand{\thefootnote}{$\star$}
\footnote{
          $\mu_n(\cdot)$ corresponds to the trend function $f$ in Wormald's DEM \cite{Wormald Ann. App Prob}, in which it is assumed Lipschitz continuous.\vspace*{.5em} 
         } 
     \addtocounter{footnote}{-1}     
\endgroup%
\!and $\nu_n (u_n)$ denote  $\mathbb{E}(X_{n+1}^2 \mid u_n)$; they are functions of $u_n$ by Doob–Dynkin lemma.\footnote{
\vspace{0.0\baselineskip} 
Precisely, 
$
    \mathbb{E}\big( X_{n+1} \mid u_n \big) := \mathbb{E} \big( X_{n+1} \mid \sigma(u_n) \big)
$, 
$ 
    \mathbb{E}\big( X_{n+1}^2 \mid  u_n \big) := \mathbb{E}\big( X_{n+1}^2 \mid \sigma(u_n) \big).
$ 
}
$\mathbb{E} (X_{n+1} \mid u_n)$ is the directional derivative of $u_n$ along the trajectory---e.g.
\vspace{0.2\baselineskip} 
$\displaystyle\frac{ \partial u_n }{\partial z }$, 
$\displaystyle\frac{ \partial u_n }{\partial x }$   
\vspace{0.2\baselineskip} 
or along other directions (see Fig.~\ref{fig: 3-SAT} of $3$-SAT, for a heuristic illustration of $u_n$, which is concentrated around $\bar u$). 
\vspace{0.2\baselineskip}

Let $\mu_{n}'$ denote $\displaystyle \frac {d\mu_{n}(u_n)}{du_n}$
\vspace{0.2\baselineskip}
, $\mu_{n}''$ denote $\displaystyle \frac {d^2\mu_{n}(u_n)}{d^2u_n}$, $\nu_n'$ denote $\displaystyle \frac{d\nu_n(u_n)}{du_n}$, and so on. We shall use term ``smooth function" to refer that (a) the first several order derivatives (with respective to $u_n$, briefly $u$) exist and (b) they are uniformly bounded (in other words, the domain in question is bounded away from the singularities).

\begin{boldthm}[General concentration inequality]\label{con_ineq}  
If $\mathbb{E}_{i-1}\big(e^{\delta A|X_i|} \big)\le K $ and $S_n \le Cn$, then in the domain where $\mu_n$ and $\nu_n$ are smooth, the sub-Gaussian concentration inequality holds. 
\end{boldthm}
In terms of partial derivatives, the conditions in Theorem \ref{con_ineq} states, in addition to the exponential moment condition and a mild global linear-growth condition $|S_n| < Cn$, that 
\[
    \frac { \partial \mu(t,u)} {\partial u}, \quad  
    \frac { \partial^2 \mu(t,u)} {\partial u^2}, \quad 
    \frac { \partial \nu(t,u)} {\partial u}, \quad \cdots
\]     
are uniformly bounded on the domain (i.e., bounded away from the singularities). These conditions are practically easy to verify, because one can examine $\bar {u}$, the mean of $u$, which in particular satisfies a partial differential equation (PDE); e.g., Fig.~\ref{fig: 3-SAT} for 3-SAT. 
Note, the Lipschitz condition (\ref{Lipschitz continuity}) says that small changes in distant past generations won't lead to violently change in the current step.

\begin{proof}[\textnormal{\textsc{Proof} of Theorem \ref{con_ineq}}] 
By Lemma \ref{con_extension}, we only need to show that the Lipschitz condition (\ref{Lipschitz continuity}) holds, i.e.
\[
       \left| \mathbb{E}(X_n \mid u_j) - \mathbb{E}(X_n \mid u_j') \right| \le L_n \left| u_j - u_j' \right| 
                         \quad \mbox{for $j < n$   and $L_n = O(1)$}
\] 
given the smoothness of  $\mu_j$ and $\nu_j$\!\!%
\begingroup
    \renewcommand{\thefootnote}{\fnsymbol{footnote}}
    \setcounter{footnote}{1}
    \footnotemark
    \footnotetext{In most cases of interest, when $\mathbb{E} (X|u)$ is smooth, $\mathbb{E}(X^2|u)$ is smooth as well. Generally, we assume that, in domains bounded away from singularities, the law of the random variable $X$ depends smoothly on $u$, in the sense that for every smooth function $f$, $\mathbb{E}\!\left[f(X)\mid u\right]$ is smooth in $u$. 
    The regularity of $\mathbb{E}(X^2|u)$ is the price paid for dispensing with a boundedness assumption on $X$.}
\endgroup%
\!\!. In particular, when $j=0$, as at the start point of $K$-SAT formula generation, 
\[
    \Bigg|\, \mathbb{E}\left(X_n \mid u_0  =  \frac   i N \right) 
              - 
             \mathbb{E}\left(X_n \mid u_0' =  \frac{i'}{N}\right)     
    \Bigg| 
                               \le L  \frac{|i-i'|} N  
\]

We consider $\big| u_j - u_j' \big| = o(1)$ because the case where $\big| u_j - u_j' \big | = O(1) $ is trivial since $S_n \le CN$.  
We proceed by induction. The base case, i.e. 
\[ 
        \Big| \mathbb{E} (X_{j+1}|u_j) - \mathbb{E} (X_{j+1} | u_j') \Big| \le L_j \left|u_j - u_j' \right| 
        ,\quad L_{j+1} = O(1) 
\]
holds because of the smoothness of $\mu_{j}$, where the Lipschitz constant $L$ is controlled by the maximum gradient up to $j$.   

In the following, we investigate stochastic process of $\{u_n \}_{n > j}$ starting from a fixed $u_j$, which therefore treated as deterministic.  We write $\mathbb{E}_{u_j}[\, \cdot \,]$ for expectation under the conditional law given $u_j$ which is a function $u_j$. Let 
\[
      \Delta \mathbb{E}\left[ X_{j+1} \right] := \mathbb{E}_{u_j}[\, X_{j+1} \,] - \mathbb{E}_{u_j'}[\,  X_{j+1}  \,]       
\]
Taking the difference between the two trajectories gives
\[
            \Delta \bar u_{j+1}  = u_j - u_j'  +  \frac{1}{N} \Delta \mathbb{E}[\, X_{j+1} \,]
\]
Then the base estimate gives
\begin{align*}
        \left| \Delta \bar u_{j+1} \right| & \le \left| u_j - u_j' \right|  +  \frac{1}{N} L_j \left| u_j - u_j' \right|  
                               \\[3pt]
                               & = \left( 1 + \frac{L_j}{N}  \right) \left| u_j - u_j' \right | 
\end{align*}
More precisely, by the mean value theorem, there exists some $\xi_j$ between $u_j$ and $u_j'$ such that
\begin{align}\label{Delta bar u_j+1}
         \Delta \bar u_{j+1}  & =  \left( 1 + \frac{\mu'_j(\xi_j)}{N}  \right) (u_j - u_j') 
         \\[3pt]\nonumber
                              :&= L_{j+1}(u_j - u_j')
\end{align}

Bearing in mind that $\mathbb{E}(X_{n+1} \mid u_n)$ is a function of $u_n$, we consider a general case of function $f(u_n)$ which is smooth. Write $\Delta_n := u_n - \bar u_n$. We restrict $u_n$ to the neighbourhood of $\bar u_n$: 
 \[ 
      A_n = \Big \{|\Delta_n| \le  \frac{1}{N^{1/2-\varepsilon }} \Big \}, \quad 0<\varepsilon < \frac 1 2
\]   
where $\bar u$ is bounded away from the singularities and hence $\mu_n$ is also bounded away from the singularities, and $f(u_n)$ is smooth on this region. Thus, there is no concern about singularities; all relevant functions are smooth. By the inductive hypothesis, $u_n$ is concentrated around $\bar u_n$ with sub-Gaussian tail bound, implying
\[ 
        \Pr\left( |\Delta_n| > \frac{1}{N^{1/2-\varepsilon } } \right) \le  e^{-c N^\varepsilon}
\]        
which means the unrestricted and the restricted expectations differ only by an exponentially small amount $\exp\!\left(-cN^\varepsilon\right)$, and the discrepancy is negligible. This fact admits a singularity-free discussion. 

Taking expectation of the Taylor series expansion of $f_n(u_{n})$ around the mean of $u_n$ ($u_n$ itself depends on the initial $u_j$ through the evolution of the process)  gives  
\begin{equation} \label{Expanding f}  
    \mathbb{E} \left[ f(u_n) \right] = f (\bar u_n) ~ +
   \frac {1}{2}  f''(\bar u_n) \mathbb{E}\Delta_n^2
                                           + \frac{1}{3!} \mathbb{E}  \left[ f'''(\xi )\Delta_n^3 \right]
\end{equation}
Similarly, since $\mathbb{E}(X_{n+1} \mid u_n)$ is a function of $u_n$, denote by $\mu_n$, by the law of total expectation
\begin{align*}
   \mathbb{E}(X_{n+1})\, \renewcommand{\thefootnote}{$\S$}   \footnotemark  
    & = \mathbb{E} \left[ \mathbb{E}(X_{n+1} \mid u_n) \right]  
      = \,  \mathbb{E} \left[ \mu_n(u_n) \right] \\
    & = \mu_{n} (\bar u_n) ~ +
   \frac {1}{2}  \mu_{n}''(\bar u_n) \mathbb{E}\Delta_n^2 
                                           + \frac{1}{3!} \mathbb{E}  \left[ \mu_n'''(\xi )\Delta_n^3 \right]
\end{align*}
The sub-Gaussian tail bound 
yields
\begingroup
    \renewcommand{\thefootnote}{$\S$}
    \footnotetext{Given $u_j$, $\mathbb{E}(X_{n+1})$ is a function of $u_j$. We may return to this notation  $\mathbb{E} (X_{n+1} \mid u_j)$ later. 
                 }
\endgroup 
\[
   \mathbb{E}|\Delta_n|^p \le \frac {C_p} {N^{p/2}}, \;   \mbox {for  $p \ge 1$ } 
                    \text{
                          (see Appendix \ref{sub-Gaussian<=> moment ≤ C_p/N^{p/2}} for a proof) 
                         } 
\]
By the hypothesis of ``smoothness'', $\mu_{n}'''(\xi)$ exists. Hence
\[
     \Big|\, \mathbb{E} \left[ \mu_{n}'''(\xi )\Delta_n^3 \right] \Big|   \lesssim  \;    {N^{-3/2}}
\]
and therefore  
\begin{equation} \label{E(X_{n+1})}
  \mathbb{E} (X_{n+1}) = \mu_{n} (\bar u_n) ~ +
               \frac {1}{2} \mu_{n}''(\bar u_n) \mathrm{Var}_n   +  \sim \frac{1}{N^{3/2}}
\end{equation}
where $\mathrm{Var}_n$ denotes $\mathbb{E} (u_n - \bar u_n)^2$.
Similarly, for $\nu$ we have
\begin{align}\label{EX{n+1}^2)}
  \mathbb{E} (X_{n+1}^2)   =& \nu_{n} (\bar u_{n})   + \frac {1}{2}   \nu_{n}''(\bar u_{n}) \mathrm{Var}_n
               +   \sim \frac{1}{N^{3/2}}
\end{align}
To estimate $\mathrm{Var}_n$, consider
\begin{align*}
    \mathrm{Var}_{n+1} &= \mathbb{E} \Big(u+x_{j+1} \cdots  + x_{n} -  
                                                                 (u + \bar x_{j+1}+ \cdots + \bar x_{n})
         + (x_{n+1} - \bar x_{n+1}) \Big)^2
         \\
                       &:= \mathbb{E} \Big( \Delta_n 
                                + (x_{n+1} - \bar x_{n+1}) \Big)^2         
\end{align*}
which yields the following recurrence equation about $\mathrm{Var}_{n}$
\begin{equation}\label{Var_n}
    \mathrm{Var}_{n+1} =  \mathrm{Var}_n + 2\mathbb{E} [ (x_{n+1}-\bar x_{n+1})\Delta_n]
        + \frac {1}{N^2}\mathbb{E} \left(X_{n+1} - \overline X_{n+1} \right)^2
\end{equation}
Applying the law of total expectation with $\mathbb{E} [X_{n+1} \,|\, u_n] = \mu_n(u_n)$, and noting $\bar x_{n+1} = \mathbb{E}[x_{n+1}]$, the expectation in the cross term above is estimated:
\begin{align*}
   \mathbb{E} \left[ (x_{n+1}-\bar x_{n+1}) \Delta_n \right] 
   = &
   \frac{1}{N} \mathbb{E} \left[\, X_{n+1} \Delta_n\, \right] 
     - 
   \bar x_{n+1} \mathbb{E} \left[\, \Delta_n \, \right]
   \\
   = & \frac{1}{N} \mathbb{E} \Big[ \mathbb{E} [ X_{n+1} \,|\,  u_n] \, \Delta_n \Big]
   \\
   = & \frac{1}{N}  \mathbb{E} \Big[\mu_{n}(u_{n})\Delta_n \Big]
\end{align*}
By the mean value theorem, 
we have 
\[
   \mu_{n}(u_{n})   =   \mu_{n}(\bar u_n) + \mu_{n}'(\xi) \Delta_n
\] 
where $\xi$ lies between $u_n$ and $\bar u_n$. 
Let 
\[
      I_n = \left[\bar u - \frac{1}{N^{(1-\varepsilon)/2}}, \bar u + \frac{1}{N^{(1-\varepsilon)/2}} \right]
\]
\[
              m_\xi := \min_{t\in I_n} \mu'_n(t) ,  \quad M_\xi := \max_{t\in I_n} \mu'_n(t)
\]
It holds that
\[
     m_\xi  \mathbb{E}     \Delta_n^2   
              \le  \mathbb{E} \left[ \mu_{n}'(\xi) \Delta_n^2  \,\right] \le  
     M_\xi  \mathbb{E}   \Delta_n^2   
\]
Therefore, there exist $\xi_n^* \in I_n$ such that
\[
                \mathbb{E} \left[ \mu_{n}'(\xi) \Delta_n^2    \,\right]
              = \mu_{n}'(\xi_n^*) \mathbb{E} \Delta_n^2 
\]
Plugging in, 
\begin{align*}
         \mathbb{E} &\left[\,  \mu_{n}(u_{n}) \Delta_n   \, \right]    
          =   \mu_{n}'(\xi_n^*) \mathbb{E} \Delta_n^2 
\end{align*}
Hence,  
\begin{align*}
   \mathbb{E} \left[ (x_{n+1}-\bar x_{n+1}) \Delta_n  \right]     
    = & \frac{\mathrm{Var}_n}{N}   \mu_{n}'(\xi_n^*)     
\end{align*}

\vspace{0.5\baselineskip}
Now let's estimate $N^{-2} \mathbb{E} \left( X_{n+1} - \overline  X_{n+1} \right)^2 $, the third term in \eqref{Var_n}. 
Combining the identities (\ref{EX{n+1}^2)}) and (\ref{E(X_{n+1})}) and neglecting terms of order $O(N^{-3/2})$ in $\mathbb{E} \left( X_{n+1} - \overline  X_{n+1} \right)^2$, we have 
\begin{align*}
     \mathbb{E} \left( X_{n+1} - \overline  X_{n+1} \right)^2 = &  \mathbb{E} (X_{n+1}^2) -  \overline  X_{n+1}^2
    \\      
       = & \nu_{{n}} (\bar u_{n})   +
  \frac {1}{2}
   \nu_{{n}}''(\bar u_{n}) \mathrm{Var}_n -
                                 \left[ \mu_{{n}} (\bar u_{n}) + \frac {1}{2} \mu_{{n}}''(\bar u_{n}) \mathrm{Var}_n \right]^2
      \\
    = & \nu_{{n}} (\bar u_{n}) + \frac {1}{2} \nu_{{n}}''(\bar u_{n}) \mathrm{Var}_n
        -    \mu_{{n}}^2 (\bar u_{n}) - \mu_{{n}} (\bar u_{n}) \mu_{{n}}''(\bar u_{n}) \mathrm{Var}_n
   \\
    = & \nu_{{n}} (\bar u_{n}) - \mu_{{n}}^2 (\bar u_{n})
         + \left( \frac {1}{2} \nu_{{n}}''(\bar u_{n})
         - \mu_{n} (\bar u_{n})\mu_{n}''(\bar u_{n})
           \right) \mathrm{Var}_n
\end{align*}
In particular, at $j+1$ 
\[
     \mathbb{E}\left(X_{j+1} - \overline X_{j+1} \right)^2 = \mathbb{E}\left(X_{j+1}^2 \right) - \mu_{j}^2
                                      = \nu_j(u_j) - \mu_{j}^2(u_j)
                                      = \mathrm{Var}(X_{j+1})
\]
where $\mathrm{Var}(X_j)=0$ because $u_j$ (such as $u_0$) is treated as a deterministic variable.

For final variance recurrence, substituting $\mathbb{E} (X_{n+1} - \overline  X_{n+1} )^2$ and $\mathbb{E} [(x_{n+1}-\bar x_{n+1})(u_{n} - \bar u_{n})] $ back into (\ref{Var_n}), we have the recurrence relation for $n \ge j$,
\begin{align*}
    \mathrm{Var}_{n+1} &=  \mathrm{Var}_n  + \frac{2\mathrm{Var}_n  }{N} 
                                                      \mu_{n}'(\xi_n^*)    
    \\
    & 
       ~ +   
       \frac{1}{N^2}  \left[ \nu_{{n}} (\bar u_{n}) - \mu_{{n}}^2 (\bar u_{n})
         + \left( \frac {1}{2} \nu_{{n}}''(\bar u_{n})
         - \mu_{n} (\bar u_{n})\mu_{n}''(\bar u_{n})
           \right)\mathrm{Var}(u_{n})
         \right]
    \\
    = &  \mathrm{Var}_n + \frac {\mathrm{Var}_n}{N} \left( 2\mu_{n}'(\xi_n) 
       ~ +      \frac {1}{2N} \nu_{{n}}''(\bar u_{n})
              - \frac {1}{N} \mu_{n} (\bar u_{n})\mu_{n}''(\bar u_{n})              
             \right)            
   \\  & 
        + \frac{1}{N^2}  \left( 
                                  \nu_{n} (\bar u_{n}) - \mu_{n}^2 (\bar u_{n})
                          \right)
\end{align*}
and the base case is $\mathrm{Var}_j=0$. With the arguments $\bar u_{n}$ and $\xi_n^*$ omitted for brevity:
%
\begin{align*}
    \mathrm{Var}_{n+1}   = & \mathrm{Var}_{n} 
               \left(1 + \frac {2\mu_{n}'}{N} 
             \right)
      +  \frac{1}{N^2} (\nu_{n}  - \mu_{n}^2 ) 
\end{align*}
where higher order terms in $N^{-1}$ inside the parenthesis have been neglected. 
We define  
\[ 
  a_n =   1 + \frac {2\mu_{n}'}{N}  
       , \qquad   b_n = \frac{\nu_{n}  - \mu_{n}^2 }{N^2} 
\]
Then, the variance recurrence relation for $n \ge j$ takes the form 
\[
       \mathrm{Var}_{n+1}   = a_n   \mathrm{Var}_{n} + b_n
\] 
It follows  
\begin{align*}
       \mathrm{Var}_{n+1}
                = & \mathrm{Var}_{j} \prod_{k=j}^{n} a_k + \sum_{r=j}^{n} \left( b_r \prod_{k=r+1}^{n} a_k \right)
                \qquad (\text{where $\mathrm{Var} _j=0$})
            \\
                = & \frac {1}{N^2} \left[  \sum_{r=j}^{n}    \left( \nu_{r}  - \mu_{r}^2 \right)
                          \prod_{k=r+1}^{n}  \left(1 + \frac {2\mu_{k}'}{N} \right)  
                                   \right]                                     
\end{align*}

\vspace{0.2\baselineskip} 
From the above expression, it is easy to check that, on any domain where $\mu_r$, $\nu_r$ and $\mu'_r$ are uniformly bounded for
\vspace{0.15\baselineskip}
$r=1, 2, \ldots, n$, $\bm{\mathrm{Var}_n = O(N^{-1})}$. Heuristically, since each step contributes $O(N^{-2} )$ to the total variance, provided each variance exists,   $n=O (N )$ steps result in a cumulative variance of $O ( N^{-1}  )$. 

We now estimate the derivative of $\mathrm{Var}_n$ and show that $\bm { \mathrm{Var}_n' } $ \textbf{is also} $\bm {O(N^{-1}) }$, which is needed for the rest of the proof and will also be used later in deriving the $K$-SAT PDE. Employing the expression of $\mathrm{Var}_{n+1}$, hereafter denoted by $V_{n+1}$, we have
\[
       V_{n+1} - V_n = O(1) \frac {1} {N^2}
\]
More generally, over $\Delta n$ steps, 
\[
       V_{n+\Delta n } - V_n = O\left( \frac {\Delta n} {N^2} \right) 
\]
and thus
\[
           \frac{V_{n+\Delta n}  - V_n}{\Delta y} = \frac {O(1)} N, \qquad \Delta y = \frac{\Delta n}{N}
\]
and thus
\[
            \frac{dV}{dy} = O(N^{-1})
\]
Note that in the preceding argument, $n$ is a generic index for successive steps of the stochastic process, which therefore may be considered along either the longitude ($m$) direction, or the latitude ($i$) direction, or any diagonal path combining both directions. Therefore, in the $i$-direction (corresponding to an increase in the number of fixed variables, in $K$-SAT), and in the $m$-direction (corresponding to the addition $K$-clauses), respectively, we have
\begin{equation}\label{Var'}
    \frac  {\partial V_n} { \partial x } =  O(N^{-1}), \qquad \frac  {\partial V_n} { \partial z } =  O(N^{-1})
\end{equation}
where
\[ 
      \Delta x = \frac{\Delta i} N, \qquad   \Delta z = \frac{\Delta m} N
\] 
Let $\Lambda_j :=  \displaystyle \frac {d \bar u_j} {d y}$ (in the $K$-SAT setting, ${d \bar u_0}/{d x} = \frac 1 2$). 
By the inductive hypothesis, i.e., $\Delta \bar u_n = L_n (u_j - u_j') $, the chain rule calculation gives
\[    
      \frac {dV_n}{d\bar u_n} 
      = \frac{1}{\Lambda_j}\frac {1}{L_n}  \frac{\partial V_n }{\partial y}      
      = O\left(\frac{1}{N}\right)
\]

Consequently, by (\ref{Expanding f}) we have
\begin{align*}  
  \Delta \mathbb{E} [f(u_n) ]  &=  \Delta f(\bar u_n) ~   +  ~  \Delta  \frac {1}{2} f_n'' \mathrm{Var}_n 
                                                                                                 ~ + ~  \frac{C}{N^{3/2}}\,  
       \\                                                                                                  
        & =   f_n' \Delta \bar u_n
           +   \frac 1 2 f_n''' \mathrm{Var}_n \Delta \bar u_n  +
                    \frac {1}{2} f_n''  \mathrm{Var}_n' \Delta \bar u_n ~ + ~  O(\Delta \bar u_n^2)
                    ~ + ~  \frac{C}{N^{3/2}}\, 
\tag{Expanding}                    
        \\
               & =   f_n' \Delta \bar u_n  
                         + \frac {1}{2} f_n''' O \Big(\frac {1}{N} \Big) \Delta \bar u_n  +
                        \frac {1}{2}  f_n''O\Big(\frac {1}{N} \Big) \Delta \bar u_n   ~ +  ~ O(\Delta \bar u_n^2)
                        ~ + ~  \frac{C}{N^{3/2}}\, 
        \\[5pt]
        & =   \Delta \bar u_n \Big( f_n' +  O\left( N^{-1} \right)  + O(\Delta \bar u_n)  \Big)
                  ~ + ~  \frac{C}{N^{3/2}}\, 
\end{align*} 
where $f_n'$ is shorthand for $f_n'(\xi)$, with $\xi$ lying between $\bar u_n$ and $\bar u_n'$.

Similarly, from (\ref{E(X_{n+1})}) where $f(\cdot) = \mu(\cdot)$,   
\[
  \Delta \mathbb{E} (X_{n+1})   =   \Delta \mu_{n} ~   +  ~  \Delta  \frac {1}{2} \mu_{n}'' \mathrm{Var}_{n} 
                                                                                                 ~ + ~  \frac{C}{N^{3/2}}\, . 
\] 
yielding
\begin{equation}\label{Delta E X_n+1}
     \Delta \mathbb{E} (X_{n+1})  
     =   \Delta \bar u_{n} \Big(\mu_{n}' + O\left(N^{-1}\right) + O\left(\Delta \bar u_{n}\right) \Big) 
     ~ + ~  \frac{C}{N^{3/2}}\, 
\end{equation}
It follows from  
\[
  u_{n+1} = u_{n} + \frac{X_n}{N}, \qquad    \bar u_{n+1}= \bar u_{n} +  \frac{1}{N}  \mathbb{E}(X_{n+1})
\] 
that
\begin{align*}
 \Delta \bar u_{n+1}&= \Delta \bar u_{n} + \frac {1}{N}  \Delta \mathbb{E}(X_{n+1}) \\
                  &= \Delta \bar u_{n} + 
                          \frac {1}{N}  \Delta \bar u_{n} \Big(\mu_{n}' + O(N^{-1}) + O(\Delta \bar u_{n}) \Big) 
                           + O\left(N^{-5/2}\right)
              \\
                  &= \Delta \bar u_{n} \left(  1 + \frac {1}{N}   \mu_{n}' + \frac{o(1)}{N}   \right) 
                   +  O\left(N^{-5/2}\right)
\end{align*}   
By iteration on $n$ and neglecting the term of $o(1)N^{-1}$, it follows
\begin{align*}
           \Delta \bar u_{n+1} &= \Delta \bar u_{j+1} \prod_{k=j+1}^n \left(1 + \frac{\mu'_k}{N} \right)  
            \\[3pt] 
\tag{By \eqref{Delta bar u_j+1}}    &= (u_j - u_j') L_{j+1}  \prod_{k=j+1}^n \left(1 + \frac{\mu'_k}{N} \right)  
            \\[3pt]                                               
                                    :&= L_{n+1} (u_j - u_j')
\end{align*} 
where the additive term $O\big(N^{-3/2}\big)$ is omitted. 
Incorporating the induction base, we have
\[
           \Delta \bar u_{n+1} \le \prod_{k=j}^n \left(1 + \frac{\mu'_k}{N} \right) \cdot  \left| u_j - u_j' \right|  
\] 
In particular, when $\left| u_0 - u_0' \right| = O\left(N^{-1}\right) $
\[
           \Delta \bar u_{n+1} = O\left(N^{-1}\right)
\] 
This is consistent with the fact that the solution of a PDE is continuous with respect to its initial condition (see, e.g.,  Fig.~\ref{fig: 3-SAT}). 

Finally, writing
\[
            M_n :=   \sup_{u}|\mu_n'(u)|
\]
and noting 
\begin{align*}
  \Big| \mathbb{E} (X_{n+1}  \mid u_j) - \mathbb{E} (X_{n+1}  \mid u_j') \Big|
                &=   \Big|   
                              \Delta   \mathbb{E}  \left[ \mu_n(u_n)   \right]  
                    \Big|
               \\[4pt]     
                &= \big|  \Delta \bar u_n \left( \mu_n'(\xi) + o(1) \right)  \big |   
                                                 \quad (\text{see \eqref{Delta E X_n+1}})
                \\[4pt]
                &\le
                   \big(M_n + o(1) \big) \big| \Delta \bar u_n  \big | 
\end{align*}
we obtain,
\begin{align*}   
                \Big|\,   \mathbb{E} (X_{n+1}  \mid u_j) - \mathbb{E} (X_{n+1}  \mid u_j') 
               \,\Big| 
            & \le 
                \big(M_n + o(1) \big) L_{n} \left| u_j - u_j' \right|  
\end{align*}
with  $L_n = O\left(L_{j+1}\right)$ when $n = O(N)$.  
This concludes the proof.
\end{proof}

\section{\textbf Concentration inequalities for BRW}\label{Concentration inequalities for BRW}

\subsection{Branching random walk (BRW)}\label{Branching random walk (BRW)}
In the literature, a branching random walk on the real line is described as follows (see, for example, \cite{Biggins_fst_lst} \cite{Biggins_Mtgl} \cite{Biggins_Chernoff}). In generation zero, an initial particle is located at the origin on the real line $\mathbf{R}$. It splits into a random number of child particles who form generation one. The children's displacements, relative to their parent, correspond to a point process on $\mathbf{R}$. The children in turn split too to form the second generation, and so on. If the average split number (branching factor) is greater than one, then with positive probability the number of the descendants grows exponentially through generations. Classical BRW studies typically address models where offspring and their displacements are independent across generations (e.g. \cite{Biggins_Mtgl}, \cite{Biggins_Chernoff} \cite{BIGGINS_et1997}, \cite{Hu&Shi2009}, \cite{Gao_BRW2014}). The law of large numbers and central limit theorem type results about the distribution of position are established, under some conditions of independence (for instance, i.i.d. branching and walking).

Considering i.i.d. offspring (and hence independent of the parent's position), and i.i.d. displacement, Harris \cite{Harris63} conjectured that the distribution of the descendants' position in the $n$th generation approaches a Gaussian distribution. This conjecture was proved by \cite{Stam66}, \cite{Kaplan_II}. Its extended generation-dependent versions, where offspring and displacement distribution depend on generation $n$, were proved by \cite{Biggins_CLT}, \cite{Klebaner82}, to mention a few. Problems where offspring and displacement are dependent on their parents' positions are studied through model-specific approaches though (e.g. \cite{Yoshida}). Concerning the deviation from the expectation, \cite{Biggins_Chernoff} showed that for any fixed $\delta > 0$, the proportion of particles whose scaled position is less than $\mu - \delta$ (or greater than $\mu + \delta$) converges almost surely to zero. In another direction of quantifying the rarity of the large deviation, extremum has been well studied (e.g. \cite{DEKKING1991} \cite{MCDIARMID_95}, \cite{BACHMANN2000}, \cite{BRAMSON&Zeitouni 2009} \cite{Hu&Shi2009}, \cite{Addario2009}  \cite{Aidekon 2013} ). \cite{DEKKING1991} showed tightness for $M_n - \mathbb{E} M_n$, where $M_n$ is the minimal position in $n$th generation. \cite{MCDIARMID_95} gave probability bound for the deviation of $M_n$, i.e. $\Pr(M_n - \mathrm{Med}_{n} > x) < e^{-\delta x}$, where $\mathrm{Med}_n$ is the median of $M_n$. A Chernoff bound for deviations of the descendants' positions in the $n$th generation has only recently been established for the aforementioned i.i.d.\ case \cite{C. Liu Chernoff Bound BRW}. To our knowledge, this is the only existing work on Chernoff type concentration inequalities for BRW.

\subsection {Motivation and notation}\label{Motivation and notation}
Concentration inequalities estimate tail probabilities:
\begin{equation} \label{RW Pr}
    \int_{|S_n - \mathbb{E}(S_n)| > \lambda}  p_{1}(x_1 | u_0)p_{2}(x_2 | u_1) \cdots p_n(x_n| u_{n-1})dx_1dx_2 \cdots dx_n
\end{equation}
where $S_n = X_1 + X_2 + \cdots + X_n$, $p_{1}(x_1 | u_0)p_{2}(x_2 | u_1) \cdots p_n(x_n| u_{n-1})$ is the probability density of $(x_1, \ldots, x_n)$ (here and throughout we use $p_i(x_i)$ for $p_i(x_i| u_{i-1})$ for notational simplicity).

In the algorithm analysis and the study of random graphs, which are modelled as stochastic processes of graph evolution (see e.g. \cite{Karp_Sipser 1981}, \cite{Wormald1}, \cite{Wormald2}, \cite{Wormald Ann. App Prob}, \cite{Wormald in Warsaw} and references therein), the differential equation method (DEM) originated in the 1980s and generalized by in \cite{Wormald Ann. App Prob}. Recently DEM has obtained some important results and received increasing attention. The procedure of DEM (see \cite{Wormald1} for a general framework) consists of (a) verifying that the random parameters concerned concentrate around their means, and (b) finding solutions to the differential equations satisfied by the means.
 
However, existing concentration inequalities, including those proved in DEM, are formulated in the setting of random-walk processes (RW), whereas a large class of problems, such as the sequential generation of ensembles of $K$-SAT and $q$-COL instances, are modelled as branching random walks (BRWs). 

Let's first fix some notation. It is suggested \cite{C. Liu Chernoff Bound BRW} that a BRW be defined by a sequence of pairs
\[
            \Big( a_i(u_{i-1}), p(x_i \mid u_{i-1}) \Big)_{i=1, 2, \hdots } 
\]
where the branching factor $a_i$ (also called birth-rate) is the expectation of offsprings which depends on the birthplace $u_{i-1}$, and
 \[
      p(x_i \mid u_{i-1})dx : = \Pr(X_i \in dx \mid \sigma(u_{i-1}))
 \]
is the probability density function (pdf) of children's displacements which depends on the birthplace as well. Here $a_i$ can be either greater than 1 (so called supercritical) or less than 1 (subcritical).
The probability law of the spine (a single lineage of descendants) is (\cite{C. Liu Chernoff Bound BRW})
\begin{equation}\label{law of spine} 
             \frac { \displaystyle 
                     \prod_{i=1}^{n} a_i(u_{i-1}) p_i(x_i \mid u_{i-1}) \,\renewcommand{\thefootnote}{$\star$} \footnotemark 
                   }
                   { \displaystyle 
                     \int \prod_{i=1}^{n} a_i(u_{i-1}) p_i(x_i \mid u_{i-1}) \, dx_1 dx_2 \cdots dx_n
                   } 
                   \, dx_1 dx_2 \cdots dx_n 
\end{equation}
which can be interpreted as the proportion of spines $(x_1, x_2, ..., x_n)$ in the whole forest. In the context of $K$-SAT, it is this tail probability that we are concerned with
\begingroup
\renewcommand{\thefootnote}{$\star$}
  \footnotetext{$\displaystyle u_{i-1} = x_0 + x_1 + \cdots + x_{i-1}$. 
               }
  \addtocounter{footnote}{-1}%
\endgroup%
\begin{equation} \label{BRW Pr}
    \int_{\left|\, S_n - \mathbb{E}(S_n) \, \right| > \lambda} 
             \frac { \displaystyle 
                     \prod_{i=1}^{n} a_i(u_{i-1}) p_i(x_i \mid u_{i-1}) 
                   }
                   { \displaystyle 
                     \int \prod_{i=1}^{n} a_i(u_{i-1}) p_i(x_i \mid u_{i-1}) \, dx_1 dx_2 \cdots dx_n
                   } 
                   \, dx_1 dx_2 \cdots dx_n 
\end{equation}
Except in special cases where the factor $a_i$ can be cancelled out, e.g., $a_i$ is independent of $u_{i-1}$, (\ref{BRW Pr}) cannot be reduced to the corresponding random walks' tail probability (\ref{RW Pr}). In this paper, we develop concentration inequalities directly for BRWs; this constitutes our second major contribution. Although originally motivated by the problems of $K$-SAT/$q$-COL, we believe the results are of independent interest. For the applications, the readers are referred to \cite{Liu K-SAT/q-COL}. 
   
If $a_i$ is constant, the BRW is equivalent to random walk process where $a_i~=~1$. In other words, random walk and branching random walk is mathematically equivalent if $a_i$ is constant. Probability density function (pdf) of BRW can be found in  (\ref{BRW Pr}). In $K$-SAT, for instance, $a_{n+1} = 1-u_{n}^{\!_K}$. That is, if a particle's position is $u_n$, the size of its offsprings is $1-u_n^{\!_K}$.
 
Intuitively, if the branching is ``squeezing''---i.e. if the branching factor is smaller in the regions farther from the mean trajectory ---then atypical paths receive less population weight than in the corresponding random walk. This ''squeezing'' effect makes the BRW distribution more concentration around it mean.  In other words, reducing the birth rate in distant regions pulls the population towards its mean trajectory, resulting in a more concentrated distribution. In context of $K$-SAT, the instances with fewer constraint have more descendants than those with more constraints; in particular, the reverse cannot hold.

Formally, let $Z(u_i)$ be the total progeny of the particle at position $u_i$. Given fixed $u_{i-1}$, this is a function of $X_i$, which we write as $Z(X_i)$. We make the following assumption.

\subsection{Concentration inequality for BRW}\label{Concentration inequality for BRW}
~ \\
\textbf{Assumption (negative association)}.  
If $f(\cdot)$ is a monotonically increasing function, then
\begin{equation}\label{neg_ass_ineq}
   \int f(|X_i|)\,Z(X_i)\,P(dX_i)
   \;\leq\;
   \left( \int f(|X_i|)\,P(dX_i) \right)
   \left( \int Z(X_i)\,P(dX_i) \right).
\end{equation}
If $Z(X_i)$ is decreasing in $X_i$, then the negative association inequality (\ref{neg_ass_ineq}) holds. $K$-SAT satisfies this assumption.

In a random walk (RW), a particle moves endlessly --- it never dies, never mates, and therefore never reproduces. The expectation of $(\cdot)$ is
\[
     \mathbb{E} [\, \cdot\, ]  =  \int (\, \cdot \, )P\big(dX_1 \cdots dX_n\big)
\]
where $P(dX_1 \cdots dX_n) = p_1(X_1) \cdots p_n(X_n) \ dX_1 \cdots dX_n$, with a slight abuse of notation.
In BRW, each particle produces descendants generation by generation. The population average of a quantity $(\cdot)$ is defined as
\begin{equation}\label{definition of E^M}
    \mathbb{E}^{\!^{(M)}} [\, \cdot \, ] =
    \frac { \displaystyle 
             \int (\, \cdot \,) \prod_{j=1}^M a_j \ P\big(dX_{1} \cdots dX_{\!_M}\big)
          }
          { \displaystyle
             \int \prod_{j=1}^M a_j   \ P\big(dX_{1} \cdots dX_{\!_M}\big)
          } 
\end{equation}
where $a_i = a_i(u_{i-1})$.  This is dependent on future generations, whereas in a random walk (RW) such statistics are independent of the future.

Define conditional expectation, with $P(\, \cdot \, \mid u_i) =   P\big(\, \cdot \, \mid \sigma(u_i)\big) $,
\begin{equation}\label{E_i_M}
   \mathbb{E}_{i}^{\!^{(M)}}[\, \cdot \,]    :=   \mathbb{E} {\!^{(M)}} [\, \cdot \mid u_i \,]  =
        \frac { \displaystyle 
                 \int (\, \cdot \,)   a_{i+1}a_{i+2}\cdots a_M   \ P\big(dX_{i+1} \cdots dX_{\!_M} \mid u_i \big)
              }
              { \displaystyle 
                 \int                 a_{i+1}a_{i+2}\cdots a_M   \ P\big(dX_{i+1} \cdots dX_{\!_M} \mid u_i \big)
              }
\end{equation} 
For example, 
\begin{align} \label{E^M}
\begin{split}
    &\mathbb{E}^{\,\!^{(M)}}_i [\, X_\ell  \,]  : =\mathbb{E} {\!^{(M)}} [\, X_\ell \mid u_i \,]   
    =  \frac { \displaystyle 
             \int X_\ell  \prod_{j=i+1}^M  a_j   \ P\big(dX_{i+1} \cdots dX_{\!_M} \mid u_i\big)
          }
          { \displaystyle 
             \int         \prod_{j=i+1}^M  a_j   \ P\big(dX_{i+1} \cdots dX_{\!_M} \mid u_i\big)
          }
\end{split}          
\end{align}
where $u_i = u_0 + X_1/N + \cdots + X_i/N$.  

Define
\begin{equation*}
           P^{\,\!^{(M)}}_i \Big(dX_{i+1} \cdots dX_M\Big)  =  
             \frac {   
                       a_{i+1}(u_{i}) \cdots  a_{\!_{M}}(u_{\!_{M-1}})  \ P\big(dX_{i+1} \cdots dX_M \mid  u_i \big)
                   }
                   { \displaystyle 
                  \int a_{i+1}(u_{i}) \cdots  a_{\!_{M}}(u_{\!_{M-1}})  \ P\big(dX_{i+1} \cdots dX_M \mid  u_i\big)
                   }
\end{equation*}
So the equality (\ref{E^M}) can be written as 
\begin{equation}\label{EX_l}
          \mathbb{E} {\!^{(M)}} [\, X_\ell \mid u_i \,]  
          =   \int X_\ell \  P^{\!^{(M)}}_i \Big(dX_{i+1} \cdots dX_\ell \cdots dX_M\Big) 
\end{equation}
When the branching factors $a_j$ are constant (and hence cancel out), the BRW conditional measure reduces to that of the underlying RW:  
\begin{equation*}
           P^{\!^{(M)}}_i \Big(dX_{i+1} \cdots dX_M \Big)  =  
               {   
                  \ P\big(dX_{i+1} \cdots dX_M  \mid  u_i\big)
               }              
\end{equation*}
and equality (\ref{EX_l}) becomes 
\begin{align*} 
          \mathbb{E} {\!^{(M)}} [\, X_\ell \mid u_i \,]  
         & =   \int X_\ell \  P\big(dX_{i+1} \cdots dX_M \mid  u_i\big) 
         \\
         & =   \int X_\ell \  P\big(dX_{i+1} \cdots dX_{\ell} \mid  u_i\big) 
\end{align*}
The second equality holds because, when the branching factors are constant, the future process does not affect the earlier behaviour.

By the law of iterated expectations
\[
     \mathbb{E} ^{\,\!^{(M)}} [\, \cdot \,] = \mathbb{E}^{\,\!^{(M)}} \big(\, \mathbb{E}_{i}^{\,\!^{(M)}} [\, \cdot \,] \, \big)
\]
that is, the BRW expectation satisfies the usual law of total expectation, despite being defined through weighted (branching) measures. We distinguish between two types of expectation. One is the standard (unweighted) expectation
\[ 
    \mathbb{E} \left[ \, \cdot  \mid u_i \right] := \int (\, \cdot \,) \, P(dX_{i+1} \mid u_i) 
                                                      = \int (\,\cdot\,) \, p_{i+1}(X_{i+1} \mid u_i)dX_{i+1}  
\] 
and the other is the BRW-weighted expectation $E^{\,\!^{(M)}}$ or $E^{\,\!^{(M)}}[\, \cdot \mid u_i \,]$ 

Define
\[
     a_i(u_{i-1}) := \mathbb{E} \big[ a_{i} (u)  \mid u_{i-1}  \big] 
\]
where $u_{i} = u_{i-1} + X_{i}/N $. 

The following lemma is the BRW counterpart of Lemma \ref {con_extension}. Its proof is identical and therefore omitted.

\begin{lemma}\label{con_inequality_BRW}  
 Let $A$ denote \  $1 + L \displaystyle \frac{M}{N}$,  where $M=O(N)$. If  
\begin{equation}\label{Lipschitz BRW}  
   \left| \,  
                 \mathbb{E}^{\,\!^{(M)}}_{i-1}(X_\ell \mid X_i) - \mathbb{E}^{\,\!^{(M)}}_{i-1}(X_\ell \mid X_i')  \, 
   \right | 
              \le L \displaystyle \frac{|X_i - X_i'|}{N}  
              ,\quad \text{for $ i < \ell $}                                     
\end{equation} 
and 
\begin{equation}\label{Exp moment existence}
   \mathbb{E}^{\,\!^{(M)}}_{i-1}(e^{\delta A|X_i|}) \le K \quad \textnormal{(finite MGF)}
\end{equation}  
for a constant $\delta > 0$ , $i=1, 2, \ldots  $ then for $n \le M$
\begin{equation}\label{Extension for BRW}
         \Pr[\ |S_n - \mathbb{E}^{ } (S_n) | \ge \lambda ] \le 2\exp 
               \left(
                     - \min  \left \{
                                       \frac{\delta^2}{4K^2} \frac{\lambda^2}{n}, \frac{\delta \lambda}{2} 
                            \right \}
              \right) 
\end{equation}
where $\mathbb{E}^{ } (S_n) $ is given by (\ref{definition of E^M}). 
\end{lemma}
\/*****************************************************************************************
Note, here without superscript,

\[ 
          \mathbb{E} [\, \cdot \,] := \int (\cdot)P(dX_{i+1}) = \int (\cdot)p_{i+1}(X_{i+1})dX_{i+1}  
\] 
which we call ``local'' expectation. 
We assume branching factor is a smooth function of birth-place $u$.
******************************************************************************************/
\begin{cor}\label{coro con_extension BRW} Suppose, in addition, that $|S_n| < Cn$ for some constant $C > 0$.  
%
Then, for every $\lambda > 0$,  
\begin{equation}\label{con ineq sub Gaussian BRW}
       \Pr(|S_n - \mathbb{E} S_n| \ge \lambda) \le 2\exp 
                                  \left(
                                           -  c\frac{\lambda^2}{n}  
                                  \right)                                   
\end{equation}  
where 
\[   
       c = \min \left\{\frac{\delta^2}{4K^2}, \frac{\delta}{4C}  \right \}       
\]     
\end{cor}

\textbf{Remark}.  
\begin{itemize}
\item Condition (\ref{Lipschitz BRW}) in Lemma \ref{con_inequality_BRW} is equivalent to the following Lipschitz continuity
\begin{equation*}\label{Lipschitz continuity BRW} 
       \left| \mathbb{E}^{\,\!^{(M)}} (X_\ell \mid u_i) - \mathbb{E}^{\,\!^{(M)}}(X_\ell \mid u_i') \right| \le L\left| u_i - u_i' \right| 
                         \quad \mbox{for $i < l$} 
\end{equation*}
where 
\[ 				u_i  = \displaystyle \frac{1}{N}  (X_1 + \cdots + X_{i-1} + X_i), \quad
          u_i' = \displaystyle \frac{1}{N}  (X_1 + \cdots + X_{i-1} + X_i').
\]      
 
\end{itemize}

The following theorem is the counterpart of Theorem \ref{con_ineq}, which gives easily verifiable conditions for BRW concentration inequality (\ref{con ineq sub Gaussian BRW}).

\begin{boldthm}[Concentration inequality in BRW]\label{con_ineq_BRW}   
If $\mathbb{E}_{i-1} (e^{\delta A|X_i|} )\le K $, $S_n \le Cn$, and the negative association (\ref{neg_ass_ineq}) holds, 
then, in the domain where $a_i(u)$, $\mu_{i}(u)$ and $\nu_i(u)$ are smooth, the concentration inequality holds. 
\end{boldthm}

\begin{proof}  
    See Appendix~\ref{Proof of con_ineq_BRW}. 
\end{proof}

Thus we have shown that the negative association assumption (\ref{neg_ass_ineq}) is sufficient to establish concentration inequalities for BRW, provided the underlying RW is not too scattering. Is it also necessary? Now consider the classical random walk of coin tossing, $S_n = X_1 + \cdots + X_n$ with $X_i \in [1, -1]$ (so $\mathbb{E}(S_n)$~=~0), which is also a BRW with branching factor $a_i = 1$. It is easy to check that if the branching factor is $a= 1 + \delta |u_i|$ for $\delta > 0$, then the BRW is dispersive. This suggests that, to some extent, the negative association assumption is also necessary.

In the situation of bounded $X_i$, a particle always stays around the mean trajectory as long as its birthplace $u_{i-1}$ is near that trajectory. In the unbounded case, however, a particle may enter regions near singularities (e.g., critical points in $K$--SAT) where $u_i$ is not necessarily smooth and the exponential moment condition $\mathbb{E}(e^ {\delta Xi}) < K $ may fail. Nevertheless, the good news is, since concentration inequality holds in the neighbourhood area of the mean trajectory---away from the the singularities---the outliers account for a negligible proportion of the population. Crucially, pruning these rare extremal paths does not weaken the concentration. The following theorem formalizes this observation.  

\begin{thm}\label{thm2} Let $\epsilon$ be small compared with $\delta$, $N^{1/2} \lesssim \lambda < \epsilon N$, and as before, 
\[
 S_i :=\sum_{t=1}^i X_t, \quad   \Delta S_i  :=  |S_i - \mathbb{E} S_i|, \quad  u_i := \frac{S_i}{N} 
\] 
Define the ``good neighborhood'' around the mean path
\[
        D_{i-1} := \Big( \bar u_{i-1} - \frac\lambda {N} ,  \bar u_{i-1} + \frac\lambda {N} \Big)
        , \quad \bar u_{i-1} : = \frac {\mathbb{E}S_{i-1}}{N}
\]
Assume as before
\[
        \mathbb{E} \left[ e^{\delta A|X_i|} \, \big| \, u_{i-1} \in D_{i-1} \right] \le K 
        \quad \textnormal{\textnormal{(bounded MGF)}}
\]
Then
\begin{equation}\label{around_Ex_A} 
     \begin{aligned}
      \Pr \Big( \displaystyle \max_{1 \leq i \leq n} \Big | \,S_i - \mathbb{E}[S_i] \, \Big |  \leq  \lambda    \Big)
          \geq  \prod_{i=1}^{n}\Big(1 - 2e^{-{c\lambda^2}/{i}} \Big) 
     \end{aligned} 
\end{equation}    
\end{thm} 
In our context, $N^{1/2} \ll \lambda < \epsilon N$ and $n < N$. Thus
\[
    \Pr \Big( \displaystyle \max_{1 \leq i \leq n} \Big | \,S_i - \mathbb{E}[S_i] \, \Big |  \ge  \lambda \Big) 
    \leq  2n \exp \left( -\frac{c\lambda^2}{n} \right)
\]
which is extremely (exponentially) small. 
\begin{proof}[\textnormal{\textsc{Proof} of} (\ref{around_Ex_A})]
The proof proceeds by induction. For each $i \ge 1$, set 
\[
         A_i := {\Delta S_i \le \lambda}
\]
Then, the following equivalence holds:
\[
    \displaystyle \max_{1 \leq i \leq n} \big| \,S_i - \mathbb{E}[S_i] \, \big |  \leq  \lambda 
     \Longleftrightarrow A_1  \wedge A_2 \wedge \cdots \wedge A_n 
\]
Note the implications  
\begin{align*}
    A_{i-1} &\iff |S_{i-1} - \mathbb{E}S_{i-1}| \leq \lambda \\
            &\iff \left| u_{i-1} - \bar{u}_{i-1} \right| \leq \frac{\lambda}{N} \\
            &\iff u_{i-1} \in D_{i-1}
\end{align*}

For the base case, we have $u_0 \in D_0$ where 
\[
       D_0 := \Big( u_0 - \frac\lambda {N} ,  u_0 + \frac\lambda {N} \Big);
\]
bounded MGF is satisfied for $i=1$. Thus the concentration inequality holds as 
\[
      \Pr \Big( |\, S_1 - \mathbb{E}[S_1] \,|  \leq  \lambda    \Big) 
          \geq   \Big(1 - 2e^{-{c\lambda^2}/1} \Big) 
\]
Assume the claim (\ref{around_Ex_A}) holds for all $i$ through to $n-1$, i.e,   
\[
    \Pr \Big(A_1  \wedge A_2 \wedge \cdots \wedge A_{n-1} \Big) \geq  \prod_{i=1}^{n-1}\Big(1 - 2e^{-{c\lambda^2}/{i}} \Big) 
\]
Then, given 
\[ 
      u_1 \in D_1, u_2 \in D_2, \ldots,  u_{n-1} \in D_{n-1}   
\]
the standard concentration inequality can be established for $S_n$ by the same argument as in Lemma \ref{con_extension} and Theorem \ref{con_ineq} :  
\[
    \Pr \Big( \big|\, S_n - \mathbb{E}S_n \, \big | \ge \lambda   \, \Big | \,
                                                 u_1 \in D_1, \ldots ,  u_{n-1} \in D_{n-1}  
        \Big)    
        \le 2\exp \Big( -c\frac{\lambda^2}{n} \Big)
\]
so
\[
    \Pr \Big( \big|\, S_n - \mathbb{E}S_n \, \big | \le \lambda   \, \Big | \,
                                                 u_1 \in D_1, \ldots ,  u_{n-1} \in D_{n-1}  
        \Big)     \ge    1 - 2\exp \Big( -c\frac{\lambda^2}{n} \Big)
\]
Thus 
\begin{align*}
    \Pr \Big(A_1 \wedge \cdots \wedge A_{n-1} \wedge A_n \Big) 
       &=  \Pr \Big(A_n \, \big| \, A_1 \wedge \cdots \wedge A_{n-1} \Big) \cdot   
           \Pr \Big(A_1  \wedge   \cdots \wedge A_{n-1}  \Big) 
      \\&
       \ge \Big( 1 - 2 e^ {-c\frac{\lambda^2}{n} }     \Big) \prod_{i=1}^{n-1}\Big(1 - 2e^{-{c\lambda^2}/{i}} \Big) 
\end{align*}
The claim follows. 
\end{proof}

\section {Applications}
See \cite{Liu K-SAT/q-COL} \/******************************************************************************  
                            Foc16 II.pdf  This introduces E(x,z): Pₘ⁽ⁱ⁾ = exp(N·∫_C dE(x,z)  
                             Extended Focs16 27page(II)&(I).pdf
                             ChatGPT: Solving dx/du = 0 at = 3.183 gives two fold solutions in (0,1/2):
                                                   uₛ₁ ≈ 0.2154608098       
                                                   uₛ₂ ≈ 0.2199421645 
                            For both, the corresponding x-value is essentially the same: 
                            
                                      xₛ(uₛ₁; 3.183) ≈ 0.1453232623,  
                                      xₛ(uₛ₂; 3.183) ≈ 0.1453225287
                                                   
                           Focs16 II.pdf & Full STOC 2016.pdf are two versions of K-SAT parts.
                           ********************************************************************************/
\/**********************************************************************************************************          
\clearpage
\appendix 

\section{\centering} \label{appendix:a}      
\begin{proof}[\textnormal{Proof of law of iterated expectations}]
$ \mathbb{E}^{(M)}[\, F\,] = \mathbb{E}^{(M)} \Big( \mathbb{E}^{(M)}_i[F] \Big) $

BRW joint probability measure:
\begin{equation*}
           P^{\,^{(M)}} (dX_1 \cdots dX_M)  =  
             \frac {   \prod_{k=1}^{M} m_k   \ P( dX_1 \cdots dX_{\!_M} )
                   }
                   { Z  }
\end{equation*}
\begin{equation*}
               Z : = \int \prod_{k=1}^M m_k P(dX_1 \cdots dX_M)
\end{equation*}
\end{proof}
\clearpage
***********************************************************************************************************/

\appendix 
\section{\centering}\label{sub-Gaussian<=> moment ≤ C_p/N^{p/2}}  
\begin{center}  
\boldmath
              \textbf{sub-Gaussian} $\Rightarrow$ $\mathbb{E}  |x|^p  \le \; \sim   N^{-p/2}$
\end{center} 
\begin{proof}  
By the tail-integral formula,  
\[
        \mathbb{E}  |X|^p    = \int_0^\infty  \Pr(|X|^p \ge x)   dx = \int_0^\infty  \Pr(|X| \ge x^{1/p})  dx
\]
Making the substitution of $x^{1/p} = t$, 
\[
        \mathbb{E}   |X|^p  = p \int_0^\infty  \Pr(|X| \ge t)  t^{p-1} dt
\]
By sub-Gaussian tail, 
\[
        \mathbb{E}  |X|^p \le C \cdot p \int_0^\infty  e^{-ct^2/n}   t^{p-1} dt
\]
The substitution $ t = \sqrt{n}s$ gives
\[
        \mathbb{E}   |X|^p   \le C\cdot p \int_0^\infty  e^{-cs^2}  n^{(p-1)/2} s^{p-1}  \sqrt{n}s\, ds
                                     = C \cdot p n^{p/2} \int_0^\infty  e^{-cs^2}   s^{p-1}   s\, ds
\]
Since the last integral above is finite for $p \ge 1$,  
\[
        \mathbb{E}   |X|^p     \le   C'_p   n^{p/2}  
\]
By scaling, it follows that
\[
                \mathbb{E}  \left|
                                  \frac{X}{N} 
                            \right| ^p
                           \le    {C_p}{N^{-p/2}}
\]

\end{proof}

\clearpage
\section{\centering}\label{Proof of con_ineq_BRW}  
 
\begin{thm}\label{con_ineq_BRW in Appendix}   
If $\mathbb{E}(e^{\delta A|X_i|} \mid u_{i-1})\le K $, $S_n \le Cn$, and the negative association (\ref{neg_ass_ineq}) holds, 
then, in the domain where $a_i(u)$, $\mu_{i}(u)$ and $\nu_i(u)$ are smooth, the concentration inequality holds. 
\end{thm}

\begin{proof} Let $Z_i^{(m)}(u_i)$ denote the population descended from the birthplace $u_i = u_{i-1} + \frac{X_i}N$ after $m$ generations:
\begin{align}\label{Z_i^M}
    Z_i^{(m)}(u_i) & = \int a_{i+1} \cdots a_{\!_{i+m}}  P(dX_{i+1}\cdots dX_{\!_{i+m}} \mid u_i)   ,\quad a_{i+1} = a(u_i) 
    \\ \nonumber
                   & := \int  dZ_i^{(m)}(u_i)
\end{align}
Define
\[
    \mathbb{E}_{i}^{(m)}[\, \cdot \,]    :=   \mathbb{E} {^{(m)}} [\, \cdot \mid u_i \,]   
          = \frac{\displaystyle 
                            \int (\, \cdot \,)  dZ_i^{(m)}(u_i)}{Z_i^{(m)}(u_i)}
\]
we have
\begingroup
    \renewcommand{\thefootnote}{$\ddagger$}
    \footnotetext{ $ P(dX_i \mid u_{i-1}) = p_i(X_i \mid u_{i-1})dX_i. $
                 }
    \addtocounter{footnote}{-1}%
\endgroup 
\begin{align*} 
    \mathbb{E}_{i-1}^{ (m)}[\,  e^{\delta A|X_i|}   \,] 
    &= \frac { \displaystyle 
                 \int e^{\delta A|X_i|} \,    Z_i^{(m)}(u_i) \ P(dX_i|u_{i-1})  
                                            \;  \renewcommand{\thefootnote}{$\ddagger$}   \footnotemark  
              }
              { \displaystyle 
                 \int     Z_i^{(m)}(u_i) \ P(dX_i|u_{i-1})  
              }   
\end{align*}
If $Z_i^{(m)}(u_{i-1}+X_i/N)$ is decreasing in $|X_i|$ for fixed $u_{i-1}$ (e.g., in the domain where backbone is not decreasing and so $X_i \ge 0$, in the $K$-SAT case; for $X_i < 0$, see the remark at the end of theorem proof), the negative association (\ref{neg_ass_ineq}) implies, 
\begin{align} \label{E^M < K}
    \mathbb{E}_{i-1}^{(m)}[\,  e^{\delta A|X_i|}   \,] 
    & \le \frac { \displaystyle 
                 \int e^{\delta A|X_i|}  \ P(dX_i|u_{i-1})  
                 \int   Z_i^{(m)}(u_i) \ P(dX_i|u_{i-1})  
                                            \;  
              }
              { \displaystyle 
                 \int   Z_i^{(m)}(u_i)  \ P(dX_i|u_{i-1})  
              }   
               \le K                
\end{align}
Thus, the exponential-moment condition in Lemma \ref{con_inequality_BRW}, inequality (\ref{Exp moment existence}), holds.

Now, similar to Theorem~\ref{con_ineq}, we only need to prove  
\[
       \Big |\, \mathbb{E}^{(m)}(X_\ell \mid u_i) - \mathbb{E}^{(m)}(X_\ell \mid u_i') 
                  \, \Big| \,  \le \, L|u_i - u_i'|  ~~~~ for ~ i < l
\]
i.e. $\mathbb{E}^{(m)}_{i}(X_\ell)$ is Lipschitz continuous in $u_i$.

Recall that in the proof of Theorem \ref{con_ineq}, smoothness of $\mu_{i} (u_i)$ and $\nu_i (u_i)$ was employed to establish the Lipschitz continuity of the trajectory with respect to the history. In an RW, ``local'' smoothness of $u_i$ implies the required Lipschitz continuity with respect to the history because the descendant's population is $1$. In a BRW, it remains to show that the future descendant population affects the corresponding Lipschitz constant by at most a bounded factor. In the present branching random walk (BRW) setting, the corresponding conditions are the smoothness of 
\[\mu_{i}^{(m)} (u_i) :=  \mathbb{E}_i^{ (m) }(X_{i+1}|u_i)
    , \quad  \nu_i^{ (m) } (u_i):=  \mathbb{E}_i^{ (m) }(X_{i+1}^2|u_i),
\]    
together with
\[ 
    \beta_i^{ (m) } (u_i)  := \mathbb{E}_{i}^{(m)}\! \big( a(u_{i+m}) \mid u_i \big).
\]
         
To show $\mu_{i}^{ (m) } (u_i)$, $\nu_i^{ (m) } (u_i)$ and $\beta_i^{ (m) } (u_i)$ are smooth functions of $u_i$,  we proceed by induction on the generations, recalling that the hypotheses in the theorem---the assumptions on the smoothness of $\beta_{i} = a_i(u)$, $\mu_{i}$ and $\nu_{i}$ (without the superscript ``$^{{(m)}}$'')---ensure that the base case (generation 0) holds. The following concerns generation $m$, where we use $n$ instead of $i$ for ease of comparison with the proof of Theorem~\ref{con_ineq}.

By \eqref{Z_i^M}, given birthplace $u_n$, its descendant population is
{\allowdisplaybreaks
\begin{align}\label{prod_m} 
\begin{split} 
  Z_n^{(m)}(u_{n}) &= Z_n^{(m)}(u_{n-1} + \frac{X_n}{N})  
   \\
   &= \int \Big(\prod_{i=1}^m a_{n+i} \Big) \  P(dX_{n+1} \cdots dX_{n+m} \mid u_{n}) 
              , \quad a_{n+i} = a(u_{n+i-1})   
   \\
  &=\int \Big(\prod_{i=1}^{m-1} a_{n+i}\Big) \ P(dX_{n+1} \cdots dX_{n+m-1} \mid u_{n}) \int a_{n+m} P(dX_{n+m} \mid u_{n-1+m})
   \\
   &=  \int  a_{n+m} dZ_n^{(m-1)}  
   \\    
   &=  Z_n^{(m-1)}(u_{n})    \cdot  
    \frac { \displaystyle
            \int  a_{n+m} dZ_n^{(m-1)}   
          }
          { \displaystyle
              Z_n^{(m-1)}(u_{n})   
          }  
          \\ 
    &  \hspace{4em}     \text{(one fewer branching factor in the denominator than its numerator.)}
     \\[5pt] 
   &= Z_n^{(m-1)}(u_{n})   \cdot  \mathbb{E}_n^{ (m-1) }  \left[\, a_{n+m}  \, \right]  
     \\[5pt]     
   &= Z_n^{(m-2)}(u_{n})   \cdot  \mathbb{E}_n^{ (m-2) }  \left[\, a_{n+m-1} \, \right]   
                                  \mathbb{E}_n^{ (m-1) }  \left[\, a_{n+m}   \, \right] 
     \\     
   & \vdots  
   \\   
   &= \prod_{k={0}}^{m-1} \mathbb{E}_{n}^{(k) }  \left[\, a_{n+k+1} \, \right]     ,
                                   \quad  \,\mathbb{E}_n^{(0)}[\, a_{n+1}\,] := a_{n+1} = a(u_n)  
\end{split}
\end{align}
}%
where
\[
    \mathbb{E}_{n}^{{(k)}}\left[\, a_{n+k+1} \, \right] =
    \frac {\displaystyle 
             \int  a_{n+k+1} \, dZ_n^{(k)}(u_n)
          }
          { \displaystyle 
              Z_n^{(k)}(u_n)
          }
\] 
\/*******************************************************************************************
{\allowdisplaybreaks 
\begin{align}\label{Z_{n+1}u_n for checking} 
\begin{split} 
  &Z_{n+1}(u_{n} + 0)  
  \\
   &=  \int \prod_{n+2}^m a_i \ P(dX_{n+2} \cdots dX_{\!_m} | u_{n}\!+\!0) 
   \\
   &=  \int \prod_{n+2}^{m-1} a_i \ P(dX_{n+2} \cdots dX_{\!_{m-1}} | u_{n}\!+\!0) \int a_{\!_m} P(dX_{\!_m} | u_{\!_{m-1}})
   \\
   &=  \int \prod_{n+2}^{m-1} a_i \ P(dX_{n+2} \cdots dX_{\!_{m-1}} | u_{n}\!+\!0) \, \beta(u_{\!_{m-1}})
   \\    
   &= \int \prod_{n+2}^{m-1} a_i \ P(dX_{n+2} \cdots dX_{m-1} | u_{n}\!+\!0)   \cdot  
    \frac { \displaystyle 
             \int \prod_{n+2}^{m-1} a_i \  P(dX_{n+2}  \cdots dX_{\!_{m-1}} | u_{n}\!+\!0) \beta(u_{\!_{m-1}}) 
          }
          { \displaystyle
               \int \prod_{n+2}^{m-1} a_i \ P(dX_{n+2}  \cdots  dX_{\!_{m-1}} | u_{n}\!+\!0) 
          }
     \\
   &= \int \prod_{n+2}^{m-1} a_i \ P(dX_{n+2} \cdots dX_{m-1} | u_{n}\!+\!0)  \cdot  \mathbb{E}_{n+1}^{\,\!^{(m-1)}} 
                                                          \left[\, \beta(u_{\!_{m-1}}) \mid u_{n}\!+\!0 \, \right]
     \\
   & \vdots  
   \\ 
   &= \int a_{n+2}(u_{n+2}) P(dX_{n+2} | u_{n}+0) \ \cdot \prod_{k={n+2}}^{m-1} \mathbb{E}_{n+1}^{(k)} 
                                                       \left[\, \beta(u_{k}) \mid u_{n}\!+\!0 \, \right]
   \\          
   &= \prod_{k={n+1}}^{m-1} \mathbb{E}_{n+1}^{(k) }  \left[\, \beta(u_{k} \mid u_{n}\!+\!0) \, \right]     
                                   ,
                                      \quad  \,\mathbb{E}_{n+1}^{(n+1)}[\, \beta(u_n+0)\,] = \beta(u_n + 0) =   \beta(u_n) 
\end{split}
\end{align}
}%
where
\[
    \mathbb{E}_{n+1}^{{(k)}}\left[\, \beta(u_{k}) \mid u_{n}\!+\!0)  \, \right] =
    \frac {\displaystyle 
             \int \beta(u_{k}|u_{n}\!+\!0) \, a_{n+2} \cdots a_{{k}} \ P(dX_{n+2} \cdots dX_{k} \mid u_{n}\!+\!0) 
          }
          { \displaystyle 
             \int a_{n+2} \cdots a_{{k}} \ P(dX_{n+2} \cdots dX_{k}  \mid u_{n}\!+\!0) 
          }
\] 
*****************************************************************************************/
Similarly,  
\[
    Z_{n+1}^{(m-1)}(u_{n + 1}) = \int \Big(\prod_{i=n+2}^{n+m} a_{i} \Big) \  P(dX_{n+2} \cdots dX_{n+m} \mid u_{n+1}) 
                               = \prod_{k={0}}^{m-2} \mathbb{E}_{n+1}^{(k) }  \left[\, a_{n+1+k+1} \, \right] 
\]    
Consequently, letting $X_{n+1} = 0$ yields
\begin{align}\label{Z_{n+1}u_n} 
     Z_{n+1}^{(m-1)}(u_{n}) = \prod_{k={0}}^{m-2} \mathbb{E}_{n+1}^{(k) }  \Big[ a_{n+1+k+1} \mid u_{n}\!+\!0  \Big]  
\end{align} 
where
\[
    \mathbb{E}_{n+1}^{{(k)}}\left[\, a_{n+1+k+1} \mid u_{n}\!+\!0 \, \right] =
    \frac {\displaystyle 
             \int  a_{n+1+k+1}\, dZ_{n+1}^{(k)}(u_n + 0)
          }
          { \displaystyle 
                          Z_{n+1}^{(k)}(u_n + 0)
          }
\] 
Given birthplace $u_{n+1} \left(=u_n + \frac {X_{n+1}}{N} \right)$, $\mathbb{E}_{n+1}^{(k) } \left[ a_{n+1+k+1} \right]$ is a function of $u_{n+1}$. 
Define
\begin{align} \label{mu_n^k}
\begin{split}
        \mu^{(m)}_n(u_n)  :&=  \mathbb{E}_n^{(m)} \big[\, X_{n+1}|u_n \, \big]   
            = \frac{ \displaystyle 
                   \int X_{n+1} \, dZ_n^{(m)}(u_n) 
                  }
                  { \displaystyle 
                          \int dZ_n^{(m)}(u_n) 
                  }
         \\[8pt]
         & =           
         \frac{ \displaystyle 
                \int X_{n+1} \cdot a( u_n) \, P(dX_{n+1} \mid u_n) Z_{n+1}^{(m-1)}(u_{n+1})
              }
              {    Z_{n}^{(m)}(u_{n})
              }
                , \qquad a_{n+1} = a(u_n)   
         \\
         & =  
      \frac{ \displaystyle 
              \int X_{n+1} a( u_n )  \, P(dX_{n+1} \mid u_n)    \frac { Z_{n+1}^{(m-1)}(u_{n+1}) } { Z_{n+1}^{(m-1)}(u_n + 0) }
           }
           {   Z_n^{(m)}(u_n)/Z_{n+1}^{(m-1)}(u_n + 0)
           }
\end{split}
\end{align}
\big(In random walks, $a_{n+1}a_{n+2} \cdots = 1$ and $  \mu^{(m)}_n(u_n) = \mu_{n} (u_n)$\big). 

To estimate $\mathbb{E}_n^{(m)} \big[\, X_{n+1}|u_n \, \big]$, we proceed by induction. Let 
{\allowdisplaybreaks
\begin{align*} 
        \overline X_{n+1}^{(1)}  & := \mathbb{E}_n^{(1)} \left[ X_{n+1}|u_n \right ]  
        = \frac{ \displaystyle 
                   \int X_{{n+1}}  a(u_{{n}})  \, P(dX_{{n+1}} \mid u_n)
                  }
                  { \displaystyle 
                    \int           a(u_{{n}})  \, P(dX_{{n+1}} \mid u_n) 
                  } 
        =   \int X_{{n+1}}  \, P(dX_{{n+1}} \mid u_n)   
\end{align*}
}%
which is the smooth function $\mu_n(u_n)$ (see Theorem \ref{con_ineq}). 
To estimate the population ratios in the numerator and denominator of \eqref{mu_n^k}, we look at
\[        
    \frac 
         {
            \mathbb{E}_{n}^{(k+1) }[\, \cdot \,]
         }          
         {\mathbb{E}_{n+1}^{(k) }\Big[\, \cdot \mid X_{n+1}=0 \, \Big]
         }
\] 
Recall that the average birth rate of generation $k$ after the $n\text{th}$, is
{\allowdisplaybreaks
\begin{align*}
      \mathbb{E}_{n}^{(k)} \left[\, a_{n+k+1} \, \right] & =  \frac { \displaystyle 
             \int   a_{n+k+1} dZ_n^{(k)}(u_n)
          }
          { \displaystyle
               Z_n^{(k)}(u_n)
          } 
     \\[6pt] 
     &=  \frac { \displaystyle 
      \int a_{n+1}P(dX_{n+1} | u_n) \, \int a_{n+k+1} \, dZ_{n+1}^{(k-1)}(u_{n+1})\,   
          }
          { \displaystyle  
                Z_n^{(k)}(u_n)
          }   
    \\[6pt] 
      &= \frac{  \displaystyle      \int a_{n+1}P(dX_{n+1} |u_n)  Z_{n+1}^{(k-1)}(u_{n+1})
                                                  \,  \mathbb{E}_{n+1}^{(k-1)} \left[\, a_{n+k+1} \, \right]
           }
           { 
              Z_n^{(k)}(u_n)
           } 
    \\[6pt] 
      &= \frac{  \displaystyle      \int     \,  \mathbb{E}_{n+1}^{(k-1)} \left[ a_{n+k+1} \right]  dZ_n^{(k)}(u_n)
           }
           { 
              Z_n^{(k)}(u_n)
           }
    \\
     &= \mathbb{E}_{n}^{(k) }   \left[
                                    \mathbb{E}_{n+1}^{(k-1) } \left[\, a_{n+k+1} \, \right]
                                \right] 
\end{align*} 
}%
That is, for all $n$ and $k \ge 1$
\begin{equation}\label{E_n^{(k)}}
\mathbb{E}_{n}^{(k)} \left[\, a(u_{n+k}) \, \right] =
         \mathbb{E}_{n}^{(k) }   \left[
                                    \mathbb{E}_{n+1}^{(k-1) } \left[\, a(u_{n+k}) \, \right]
                                \right] \quad 
\end{equation}
We next prove $\mathbb{E}_{n}^{(1)} \left[\, a(u_{n+1}) \, \right]$ is continuous in $u_n$, using the continuity of  $\overline X_{n+1}^{(1)}$, and then show the continuity of $\overline X_{n+1}^{(2)}$ using  $\mathbb{E}_{n}^{(1)} \left[\, a(u_{n+1}) \, \right]$, and so on: 
\[
\overline X_{n+1}^{(1)} \;\Rightarrow\; \mathbb{E}_n^{(1)} \;\Rightarrow\; \overline X_{n+1}^{(2)} \;\Rightarrow\; \mathbb{E}_n^{(2)} \;\Rightarrow\; \cdots
\] 
leading to smoothness of $\mu^{(m)}_n(u_n)$ and $\mathbb{E}_{n}^{(m)}\left[\, a(u_{n+m}) \, \right]$ for all $m \lesssim N$. The smoothness of $\nu_n^{{(k)}} (u_n)$ can be proved with smooth $\mu^{(k)}_n(u_n)$ and $\mathbb{E}_{n}^{(k-1)}\left[\, a(u_{n+k-1}) \, \right]$.

When $k=1$, for any $i$
\begin{align*}
      \mathbb{E}_i^{(1)} \left[\, a(u_i + X/N) \, \right] 
       & = \frac { \displaystyle 
                         \int a(u_i)  \, P(dX  | u_i ) \, a\Big( u_i +  \frac{X}{N} \Big) 
             }
          { \displaystyle
               \int a(u_i) \, P(dX | u_i)
          } 
       =  \int  a\Big( u_i +  \frac{X}{N} \Big)   P(dX | u_i)  
\end{align*}%
is a smooth function of $u_i$ (meaning it has continuous derivatives up to some required order) since $a(u)$ is uniformly smooth in $u$ and the conditional distribution $P(dX_{i+1} | u_i)$ varies smoothly with $u_i$. In $K$-SAT, $a(u) = 1 - u^K$. For instance, 
\begin{equation}\label{E_n^{(1)}}
    \mathbb{E}_n^{(1)} \left[\, a \Big( u_{n} +  \frac {X_{n+1}}N \Big) \, \right]    
     = a(u_n) + \frac{a'(u_n)}{N} \overline X_{n+1}^{(1)} + O(N^{-2})
\end{equation} 

Write
\begin{align*}  
      \mathbb{E}_{n+1}^{(k) } \left[\, a(u_{n+1+k})  \, \right]  := f_{n+1}^{(k)} \left (u_n + \frac {X_{n+1}}N \right)  
\end{align*}
Then,
\begin{align*}  
      \mathbb{E}_{n+1}^{(k) } \left[\, a(u_{n+1+k}) \mid X_{n+1}=0\, \right] 
                               = f_{n+1}^{(k)} \left (u_n\right)  
\end{align*}
Let $f_i^{(0)}(u_i)$ denote  $a(u_i)$ \big(namely, $\mathbb{E}_n^{(0)}[\, a(u_n)\,] = a(u_n)$\big). Write 
\[ \overline X_{n+1}^{(k)} :=  \mathbb{E}_n^{(k)} (X_{n+1}) 
\] 
thus (\ref{E_n^{(1)}}) can be written as
\begin{equation}\label{f^{(1)}}
          f_n^{(1)}(u_n) = f_{n+1}^{(0)}(u_n) + \frac{\overline X_{n+1}^{(1)}}{N} {f_{n+1}^{(0)}}'(u_n)  + O(N^{-2})
\end{equation}
and
(\ref{E_n^{(k)}}) becomes 
\[
          f_n^{(k)}(u_{n}) = \mathbb{E}_{n}^{(k)}  
                                \left[
                                       f_{n+1}^{(k-1)}(u_{n+1}) 
                                \right] 
\]

Now assume for $k \ge 2$, $f_{n+1}^{(k-1)}(u_{n+1})$ is smooth in $u_{n+1}$ for any $n$.  
Employing a Taylor expansion, we have
\begin{align} \label{E_n^k}
\begin{split}
  f_n^{(k)}(u_{n})   
  & = \mathbb{E}_n^{(k)}\!\left[ f_{n+1}^{(k-1)} \Big( u_{n} +  \frac {X_{n+1}}N \Big) \right]
  \\
  & = f_{n+1}^{(k-1)}\!\left(u_n\right)  +  \frac{\overline X_{n+1}^{(k)}}{N}{f_{n+1}^{(k-1)}}' (u_n)    +  O\!\left(N^{-2} \right)
\end{split}  
\end{align}%
which is smooth in $u_n$ for all $n$, provided $\overline{X}_{n+1}^{(k)}$ is smooth. Observe that the term $\displaystyle \frac{\overline X_{n+1}^{(k)}}{N}{f_{n+1}^{(k-1)}}' (u_n) $ involves no second order derivatives of $f_{n+2}^{(k-2)}$ which is absorbed into the $O(N^{-2})$ term.  
Thus, we have proved that for any $k$, $\mathbb{E}_{n}^{(k)} \left[\, a(u_{n+k}) \, \right]$ is smooth in $u_n$.  Hence, by (\ref{E_n^k}),
{\allowdisplaybreaks
\begin{align}\label{E/E} 
\begin{split}
   \frac {\mathbb{E}_{n}^{(k+1)} \left[\, a(u_{n+1+k}) \, \right]  }   
         {\mathbb{E}_{n+1}^{(k)}\big[\,  a(u_{n+1+k}) | X_{n+1}=0 \, \big] } 
  &=\frac{ \displaystyle 
           f_{n+1}^{(k)}\!\left(u_n \right) +  \frac{\overline X_{n+1}^{(k+1)}(u_n)}{N} {f_{n+1}^{(k)}}'\!\left(u_n\right)  
                                                                                        +  O\!\left(N^{-2} \right)
         }
         { 
           f_{n+1}^{(k)}  \left( u_n + 0 \right) 
         }    
  \\[10pt]
  & = 1 
      + \frac{\overline X_{n+1}^{(k+1)}(u_n) }{N} \cdot \frac{{f_{n+1}^{(k)}}'\!\left(u_n\right)} 
                                                          {f_{n+1}^{(k)}  \left( u_n \right) } 
      +  O\!\left(N^{-2} \right)
\end{split}
\end{align}
}%
Notice that $Z_{n+1}(u_n)$ is $Z_{n+1}\big(u_n +  {X_{n+1}}/{N}\big) \big | _{X_{n+1}=0} $\,. \,From (\ref{prod_m}) and (\ref{Z_{n+1}u_n}), we have 

(a) the population ratio in the numerator of (\ref{mu_n^k}): 
{\allowdisplaybreaks
\begin{align*}
  \frac{Z_{n+1}^{(m-1)}(u_{n+1})} {Z_{n+1}^{(m-1)}(u_{n})}  
       & = \frac{  
                  \displaystyle \prod_{k=0}^{m-2} \mathbb{E}_{n+1}^{(k)}\big[\,  a(u_{n+1+k}) \, \big]
                }
                { \displaystyle \prod_{k=0}^{m-2} \mathbb{E}_{n+1}^{(k)}\big[\,  a(u_{n+1+k}) | X_{n+1}=0 \, \big]
                } 
  \\
  & = \prod_{k=0}^{m-2}  
       \frac{ \displaystyle  
          f_{n+1}^{(k)}\!\left(u_n\right) +  {f_{n+1}^{(k)}}'\!\left(u_n\right) \frac{X_{n+1}}{N} +  O\!\left(N^{-2} \right)
            }
            { \displaystyle  f_{n+1}^{(k)}  \left( u_n  \right) 
            }
            \quad  \big(\text{smooth $f_{n+1}^{(k)}$}\big)
  \\
  & = \prod_{k=0}^{m-2}  \left(1  + \frac{X_{n+1}}{N} \frac{{f_{n+1}^{(k)}}'(u_n)}{f_{n+1}^{(k)}  \left( u_n \right) }     
                                  +   O\!\left(N^{-2} \right)
                           \right)
\end{align*}
}%
 
(b) the population ratio in the denominator of (\ref{mu_n^k}): 
{\allowdisplaybreaks
\begin{align*}
   \frac  {Z_n^{(m)}(u_n)} {Z_{n+1}^{(m-1)}(u_{n})}                 
  & = \frac{ \displaystyle \prod_{k=0}^{m-1}   \mathbb{E}_{n}^{(k)} \big[\, a(u_{n+k}) \, \big] 
               }
               { \displaystyle \prod_{k=0}^{m-2} \mathbb{E}_{n+1}^{(k)}\big[\,  a(u_{n+1+k}) | X_{n+1}=0 \, \big]
               }
   \\            
  &=\frac{ \displaystyle a(u_n) \prod_{k=0}^{m-2} \; \mathbb{E}_{n}^{(k+1)} \big[\, a(u_{n+1+k}) \, \big] 
         }
         { \displaystyle \prod_{k=0}^{m-2} \mathbb{E}_{n+1}^{(k)}\big[\,  a(u_{n+1+k}) | X_{n+1}=0 \, \big]
         }
       \\[10pt]
   &= a(u_n) \, \prod_{k=0}^{m-2}  
      \left(  1   
              +  \frac{\overline X_{n+1}^{(k+1)}(u_n) }{N} \cdot \frac{{f_{n+1}^{(k)}}'\!\left(u_n\right)} 
                                                                      {f_{n+1}^{(k)}  \left( u_n \right) } 
              +  O\!\left(N^{-2} \right)
     \right) \quad  \big(\text{by \eqref{E/E}}\big)
\end{align*}
}%
where  $  \overline X_{n+1}^{(k+1)} := \mathbb{E}_{n}^{(k+1) }  \left[X_{n+1}  \mid u_n \right]$. 

Write
\[
           A_k(u_n) =   \frac{{f_{n+1}^{(k)}}'(u_n)} {f_{n+1}^{(k)}(u_n)} \Big( X_{n+1} - \overline X_{n+1}^{(k+1)}(u_n) \Big)
\] 
Then, (\ref{mu_n^k})---that is, $\mu^{(m)}_n(u_n)$---becomes, with terms of order $O\!\left(N^{-2} \right)$ omitted,  
\begin{align*}
  \overline  X_{n+1}^{(m)}(u_n)  
      &=  \int \!  X_{n+1} \prod_{k=0}^{m-2} \! \left(
       1   
          +  
        \frac{1}{N}  A_k(u_n)
                                                 \right)   \! P(dX_{n+1}|u_n) 
\end{align*}  
where $a(u)$ cancels out.  This is a smooth function of $u_n$, since $f_{n+1}^{(k)}$, ${f_{n+1}^{{(k)}}}'$ and $\overline X_{n+1}^{(k+1)}(u_n)$ are smooth by induction for all $k+1 < m-1$,  and $f_{n+1}^{(k)}$ are bounded away from zero (recall $a(u_i) = 1 - u_i^K$ in $K$-SAT); as $m = O(N)$, expressions of the form $\exp\!\left(N^{-1}\sum  A_k(u_n) \right)$ are uniformly bounded,  
In the same way, we can prove the smoothness of $\nu_n^{{(m)}} (u_n)$ and $\mathbb{E}_{n}^{(m)} \left[\, a(u_{n+m}) \, \right]$. 

\textbf{Remark}. 
\begin{itemize} 
\item   In the $K$-SAT case, in the domain of 1-clause decreasing (i.e., $X_i < 0$) and $m = o(N)$ 
\begin{align*}
\mathbb{E}_{i-1}^{(m)}[\,  e^{\delta A|X_i|}   \,]  
              &= \int \!  e^{\delta A|X_i|}  \prod_{k=0}^{m-2} \! \left(
                   1   
                      +  
                    \frac{1}{N}  A_k(u_{i-1})
                                                             \right)   \! P(dX_{i}|u_{i-1})  
              \\
              &
                \le K + o(1)
\end{align*}

Here the $o(1)$ term is uniform. Thus on a vertical segment of the process where only 1-clause decreases occur, with $m = o(N)$, the exponential moment is bounded uniformly along the segment, and the concentration inequality holds. Immediately after passing this segment, the increase in the number of $K$-clauses preserves the negative association condition. 

\hspace{1em} On the other hand, the effect of decreasing a 1-clause can be offset (the case where $X_i < 0$), through negative association, by adding several 3-clauses\,\footnote{Roughly, by $k$ 3-clauses such that 
  $(1 + u)(1 - u^3)^k \le 1$
}, and afterward, the negative association condition is preserved and strengthened through the addition of 3-clauses. Over the entire 2D area of $N \times N$ points, such 1-clauses account for only $O(N)$, since the total number 1-clauses is $O(N)$. Therefore, the finite  exponential moment condition holds on a proportion $1 - O(N^{-1})$ of the whole domain. 
 
\end{itemize}%

\noindent The rest of the proof is the same as in Theorem \ref{con_ineq}. It is worthy noting that to complete the proof we need to use
\[
   \mathbb{E}^{ (m) } |u_n - \bar u_n|^p \le \frac {C_p} {N^{p/2}} ~~~~ \mbox { for  $p \ge 1 $
                                           }
\]
This in turn requires concentration inequality for $x_1 + x_2 + \cdots  x_n$ which holds by induction. 
\end{proof}

\end{document}